\documentclass[aos,preprint]{imsart}

\setattribute{journal}{name}{}

\usepackage[square,numbers]{natbib}
\RequirePackage[colorlinks,citecolor=blue,urlcolor=blue]{hyperref}
\usepackage{lipsum}                   
\usepackage{xargs}                    
\usepackage[pdftex,dvipsnames]{xcolor}
\usepackage{mathrsfs,amsthm,amsmath}
\usepackage[colorinlistoftodos,prependcaption,textsize=tiny]{todonotes}
\newcommandx{\unsure}[2][1=]{\todo[linecolor=red,backgroundcolor=red!25,bordercolor=red,#1]{#2}}
\newcommandx{\change}[2][1=]{\todo[linecolor=blue,backgroundcolor=blue!25,bordercolor=blue,#1]{#2}}
\newcommandx{\info}[2][1=]{\todo[linecolor=OliveGreen,backgroundcolor=OliveGreen!25,bordercolor=OliveGreen,#1]{#2}}
\newcommandx{\improvement}[2][1=]{\todo[linecolor=Plum,backgroundcolor=Plum!25,bordercolor=Plum,#1]{#2}}
\newcommandx{\thiswillnotshow}[2][1=]{\todo[disable,#1]{#2}}

\usepackage{amssymb,multirow}
\usepackage{amsfonts}
\usepackage{stmaryrd,color}
\usepackage{dsfont}
\usepackage[scriptsize]{subfigure}
\usepackage{ucs}
\usepackage{mathtools}
\usepackage[utf8x]{inputenc}
\usepackage{amsmath}
\usepackage{amsfonts}
\usepackage{amssymb}
\usepackage{amsthm}
\usepackage{fontenc}
\usepackage{graphicx,color}

\usepackage{bbm}
\usepackage{bm}
\RequirePackage[OT1]{fontenc}

\startlocaldefs

\newcommand{\D}{{\rm d}}

\newcommand{\bbC}{\mathds{C}}

\newcommand{\Z}{\mathds{Z}}

\newcommand{\R}{\mathds{R}}
\newcommand{\E}{\mathds{E}}
\newcommand{\N}{\mathds{N}}
\renewcommand{\P}{\mathds{P}}

\newcommand{\Ecr}{\mathscr{E}}

\newcommand{\Dcr}{\mathscr{D}}

\newcommand{\Pcr}{\mathscr{P}}

\newcommand{\Ncr}{\mathscr{N}}

\newcommand{\ind}{\mathds{1}_}
\newcommand{\Var}{{\rm Var}}

\newcommand{\MSE}{{\rm MSE}}

\newcommand{\TV}{\mathsf{TV}}

\newcommand{\pol}{\mathsf{P}}
\newcommand{\tpol}{\mathsf{TP}}
\newcommand{\cont}{\mathsf{C}}
\newcommand{\conteven}{\mathsf{CE}}

\theoremstyle{plain}
\newtheorem{thm}{\textsc{Theorem}}

\newtheorem{lem}{\textsc{Lemma}}
\newtheorem{cor}{\textsc{Corollary}}

\newtheorem{prp}{\textsc{Proposition}}

\allowdisplaybreaks

\endlocaldefs

\begin{document}

\begin{frontmatter}
\title{Optimal disclosure risk assessment}
\runtitle{Optimal disclosure risk assessment}

\begin{aug}
\author{\fnms{F.} \snm{Camerlenghi}\thanksref{m1,m4}\ead[label=e1]{federico.camerlenghi@unimib.it}},
\author{\fnms{S.} \snm{Favaro}\thanksref{m2,m4}\ead[label=e2]{stefano.favaro@unito.it}},
\author{\fnms{Z.} \snm{Naulet}\thanksref{m5}\ead[label=e4]{znaulet@utstat.toronto.edu}}
\and
\author{\fnms{F.} \snm{Panero}\thanksref{m3}\ead[label=e3]{francesca.panero@stats.ox.ac.uk}}

\thankstext{m4}{Also affiliated to Collegio Carlo Alberto, Torino, Italy}
\runauthor{Camerlenghi, Favaro, Naulet and Panero}

\affiliation{University of Milano - Bicocca\thanksmark{m1}, University of Torino\thanksmark{m2}, University of Toronto\thanksmark{m5} and University of Oxford\thanksmark{m3}}

\address{Department of Economics, Management and Statistics\\
Milano, 20126\\
Italy\\
\printead{e1}}

\address{Department of Economics and Statistics\\
Torino, 10134\\
Italy\\
\printead{e2}}

\address{Department of Statistical Sciences, Sidney Smith Hall\\ 
Toronto, M5S 3G3\\
Canada\\
\printead{e4}}

\address{Department of Statistics\\
Oxford, OX1 3LB,\\
United Kingdom\\
\printead{e3}}
\end{aug}

\begin{abstract}
Protection against disclosure is a legal and ethical obligation for agencies releasing microdata files for public use. Consider a microdata sample of size $n$ from a finite population of size $\bar{n}=n+\lambda n$, with $\lambda>0$, such that each record contains two disjoint types of information: identifying categorical information and sensitive information. Any decision about releasing data  is supported by the estimation of measures of disclosure risk, which are functionals of the number of sample records with a unique combination of values of identifying variables. The most common measure is arguably the number $\tau_{1}$ of sample unique records that are population uniques. In this paper, we first study nonparametric estimation of $\tau_{1}$ under the Poisson abundance model for sample records. We introduce a class of linear estimators of $\tau_{1}$ that are simple, computationally efficient and scalable to massive datasets, and we give uniform theoretical guarantees for them. In particular, we show that they provably estimate $\tau_{1}$ all of the way up to the sampling fraction $(\lambda+1)^{-1}\propto (\log n)^{-1}$, with vanishing normalized mean-square error (NMSE) for large $n$. We then establish a lower bound for the minimax NMSE for the estimation of $\tau_{1}$, which allows us to show that: i) $(\lambda+1)^{-1}\propto (\log n)^{-1}$ is the smallest possible sampling fraction; ii) estimators' NMSE is near optimal, in the sense of matching the minimax lower bound, for large $n$. This is the main result of our paper, and it provides a precise answer to an open question about the feasibility of nonparametric estimation of $\tau_{1}$ under the Poisson abundance model and for a sampling fraction $(\lambda+1)^{-1}<1/2$.
\end{abstract}

\begin{keyword}[class=MSC]
\kwd[Primary ]{62G05 }
\kwd[; Secondary ]{62C20}
\end{keyword}

\begin{keyword}
\kwd{disclosure risk assessment} \kwd{microdata sample} \kwd{nonparametric inference} \kwd{optimal minimax procedure} \kwd{Poisson abundance model} \kwd{polynomial approximation}
\end{keyword}

\end{frontmatter}

\section{Introduction}

Protection against disclosure is a legal and ethical obligation for agencies releasing microdata files for public use. Any decision about release requires a careful assessment of the risk of disclosure, which is supported by the estimation of measures of disclosure risk (e.g., \citet{Wil(01)}). Consider a microdata sample $\bm{X}(n)=(X_{1},\ldots,X_{n})$ from a finite population $\bm{X}$ of size $\bar{n}>n$, and without loss of generality assume that each sample record $X_{i}$ contains two disjoint types of information for the $i$-th individual: identifying information and sensitive information. Identifying information consists of the values of a set of categorical variables, which might be matchable to known units of the population. A risk of disclosure arises from the possibility that an intruder might succeed in identifying a microdata unit through such a matching and hence be able to disclose the sensitive information on this unit. To quantify the risk of disclosure, microdata sample records are cross-classified according to potentially identifying variables, i.e., $\bm{X}(n)$ is partitioned in $K_{n}\leq n$ cells with corresponding frequency counts $(Y_{1}(\bm{X},n),\ldots,Y_{K_n}(\bm{X},n))$ such that $\sum_{1\leq i\leq K_{n}}Y_j (\bm{X},n)=n$, where 
 $Y_j (\bm{X},n)$ denotes the frequency of the $j$-th cell out of the sample $\bm{X}(n)$. A risk of disclosure arises from cells in which both sample frequencies and population frequencies are small. Of special interest are cells with frequency $1$ (singletons or uniques) since, assuming no errors in the matching process or data sources, for these cells the match is guaranteed to be correct. This has motivated inference on measures of disclosure risk that are functionals of the number of singletons, the most common being the number $\tau_{1}$ of sample singletons which are also population singletons. See, e.g., \citet{Bet(90)} and \citet{Ski(94)} for a thorough discussion on measures of disclosure risk. 

The Poisson abundance model is arguably the most natural, and weak, modeling assumption to infer $\tau_{1}$. If $\bar{n}=n+\lambda n$, with $\lambda>0$, it assumes that: i) the population records $(X_{1},\ldots,X_{n+\lambda n})$ can be ideally extended to a sequence $\bm{X}:= (X_{i})_{i\geq1}$, of which $\bm{X}(n)$ is an observable subsample; ii) the $X_i$'s are independent and identically distributed according to an unknown distribution $(p_j)_{j  \geq 1}$, where $p_{j}$ is the probability of the $j$-th cell in which $\bm{X}$ may be cross-classified; iii) the sample size is a Poisson random variable $N$ with mean $n$, in symbols $N \sim {\rm Poiss} (n)$. Then sample records $\bm{X}(N)=(X_1,\ldots , X_N)$ result in  $K_{N}$ cells with frequencies $(Y_1(\bm{X},N), \ldots , Y_{K_N}(\bm{X},N))$ such that $Y_j(\bm{X},N)\sim {\rm Poiss} (np_{j})$ for $j=1,\ldots,K_{N}$, $Y_{j_1}(\bm{X},N)$ is independent of $Y_{j_2}(\bm{X},N)$
for any $j_1 \not = j_2$, and $\sum_{1\leq j\leq K_N}
Y_j (\bm{X}, N) = N$. As discussed in Section 2.4 of \citet{Ski(02)}, nonparametric estimation of $\tau_{1}$ under the Poisson abundance model is an intrinsically difficult problem. It shares the well-known difficulties of the classical problem of estimating the number of unseen species (e.g., \citet{Goo(56)}, \citet{Efr(76)}, \citet{Orl(17)}). In particular, nonparametric estimators of $\tau_{1}$ may be ``very unreasonable'' since they are subject to serious upward bias and high variance for small sampling fractions of the population, i.e. for $(\lambda+1)^{-1}<1/2$. To overcome these issues, in the last three decades stronger modeling assumptions have been considered. These studies resulted in a range of parametric and semiparametric approaches, both frequentist and Bayesian, to infer $\tau_{1}$, e.g.,  \citet{Bet(90)}, \citet{Sam(98)}, \citet{Ski(02)}, \citet{Rei(05)}, \citet{Rin(06)}, \citet{Ski(08)}, \citet{Man(12)}, \citet{Man(14)}, \citet{Car(15)} and \citet{Car(18)}. 

In this paper, we first study nonparametric estimation of $\tau_{1}$ under the Poisson abundance model for sample records. Given a collection of sample records $(X_{1},\ldots,X_{n})$ from
the population $(X_{1},\ldots,X_{n+\lambda n})$, we introduce a class of
nonparametric linear estimators of $\tau_{1}$ that are simple, computationally
efficient and scalable to massive datasets. We show that our estimators admit an
interpretation as (smoothed) nonparametric empirical Bayes estimators in the
sense of \citet{Rob(56)}, and we prove theoretical guarantees for them that hold
uniformly for any distribution $(p_{j})_{j\geq1}$. In particular, we show that
the proposed estimators provably estimate $\tau_{1}$ all of the way up to the sampling fraction 
$(\lambda+1)^{-1}\propto (\log n)^{-1}$, with vanishing normalized mean-square
error (NMSE) as $n$ becomes large. Then, by relying on recent techniques
developed in \citet{wu2015chebyshev} in the context of optimal estimation of the
support size of discrete distributions, we establish a lower bound 
for the minimax NMSE for the estimation of $\tau_{1}$. This result allows us to show
that $(\lambda+1)^{-1}\propto (\log n)^{-1}$ is the smallest possible sampling fraction of the population, and that estimators' NMSE is near optimal, in the sense of
matching the minimax lower bound, for large $n$. This is the
main result of our paper, and it provides a precise answer to the question
raised by \citet{Ski(02)} about the feasibility of nonparametric estimation of $\tau_{1}$ under the Poisson abundance model and for a sampling fraction $(\lambda+1)^{-1}<1/2$. Indeed our result shows that nonparametric estimation
of $\tau_{1}$ has uniformly provable guarantees, in terms of vanishing NMSE for large $n$, if and only if $(\lambda+1)^{-1}\propto (\log n)^{-1}$.

The paper is structured as follows. In Section 2 we introduce a class of nonparametric estimators for $\tau_{1}$, and we show that they provably estimate $\tau_{1}$ all of the way up to the sampling fraction $(\lambda+1)^{-1}\propto (\log n)^{-1}$, with vanishing NMSE as $n$ becomes large. In Section 3 we show that $(\lambda+1)^{-1}\propto (\log n)^{-1}$ is the smallest possible sampling fraction of the population, and that estimators' NMSE is near optimal for large $n$. Section 4 contains a numerical illustration of the proposed estimators. Proofs and deferred to the Appendix.

\section{A nonparametric estimator of $\tau_{1}$}
\label{sec:theory}
We consider an infinite sequence of observations $\bm{X}$, and we assume that  $\bm{X}(N)=(X_{1},\ldots,X_{N})$ is the microdata sample of random size $N$ under the Poisson abundance model. We suppose that $\bm{X}(N)$ is a subsample of  $(X_{1},\ldots,X_{M+N})$, where $M\sim {\rm Poiss} (\lambda n)$, with $\lambda >0$ and independent of $N$. In the present framework $(X_{N+1},\ldots,X_{N+M})$ may be seen as the unobservable population. When the sample records are cross-classified according to the potentially identifying variables, the sample $(X_{1},\ldots,X_{N})$ is partitioned in $K_{N}\leq N$ 
cells with corresponding frequency counts $(Y_{1}(\bm{X},N),\ldots,Y_{K_N}(\bm{X},N))$ such that $\sum_{1\leq i\leq K_{N}} Y_{j}(\bm{X},N)=N$. Hereafter we denote by $Z_i(\bm{X},N)$ the number of cells with frequency $i$, and by $Z_{\bar{i}}(\bm{X},N)$ the number of cells with frequency greater or equal than $i$, for any index $i\geq 1$. We are interest in estimating the number $\tau_{1}$ of sample uniques which are also population uniques, namely the functional
\begin{displaymath}
\tau_{1}(\mathbf{X},N,M)=\sum_{j\geq 1}\mathds{1}_{\{Y_j (\bm{X},N)=1\}}\mathds{1}_{\{Y_j (\bm{X},N+M)=1\}}.
\end{displaymath}
We recall that the frequency counts $Y_{j}(\bm{X},N)$'s are independent, and
that they are Poisson distributed with parameter $np_{j}$, where $p_{j}$ is the
unknown probability associated to the $j$-th cell, that is
$p_{j}\in[0,1]$
for $j\geq 1$ such that $\sum_{j \geq 1}p_{j}=1$. We will denote by $\bm{Y} (\bm{X},N):= (Y_1 (\bm{X},N), \ldots )$ the whole sequence of the cell's frequency count, when we are provided with a sample of size $N$.

To fix the notation, in the sequel we will write $f \lesssim g$,  for two generic functions $f$ and $g$, iff there exists a universal constant $C >0$ such that $f (x) \leq C g(x)$; we will further write $f \asymp g$ whenever both  $f \lesssim g$ and 
 $g \lesssim f$ are satisfied.
Let us denote by $\Pcr$ the set of all possible distributions over $\N$, i.e.
$\Pcr:= \{ P= \sum_{j\geq 1} p_j \delta_j :\; p_j \in [0,1],
\text{ with} \sum_{j \geq 1}p_j=1  \}$, where $\delta_j$ denotes the Dirac measure centered at $j\in \N$. 
An estimator of $\tau_1(\bm{X},N,M)$ is understood to be a measurable function  $\hat{\rho}_1 (\bm{X}(N),N)$ depending on the available sample $\bm{X}(N)$ and the actual size of the observed sample $N$. We will evaluate the performance of a generic estimator $\hat{\rho}_1 (\bm{X}(N),N)$ of $\tau_1(\bm{X},N,M)$, by its worst--case NMSE, defined as  
\begin{equation}
\label{eq:NMSE_def}
\Ecr_{\lambda , n} (\hat{\rho}_1 (\bm{X}(N),N)) := \sup_{P \in \Pcr}\frac{\E [(\hat{\rho}_1 (\bm{X}(N),N)-\tau_{1}(\bm{X},N,M))^2]}{n^2},
\end{equation}
where $\E [(\hat{\rho}_1 (\bm{X}(N),N)-\tau_{1}(\bm{X},N,M))^2]$ is the mean squared error (MSE) of $\hat{\rho}_1 $, also
denoted by $\MSE [\hat{\rho}_1  (\bm{X}(N),N)]$. The use of the NMSE \eqref{eq:NMSE_def} has been recently proposed in \citet{Orl(17)} in the context of the estimation of the number of unseen species. 

A nonparametric estimator for $\tau_1(\bm{X},N,M)$ may be deduced comparing expectations, indeed it is easy to see that:
\begin{equation}
\label{eq:comparingE}
\E [\tau_1(\bm{X},N,M)] =  \sum_{i \geq 0}(-1)^i \lambda^i (i+1) \E [Z_{i+1} (\bm{X},N)]
\end{equation}
from which we may define the following estimator 
\begin{equation}\label{eq:estimator<1}
\hat{\tau}_{1}(\bm{X}(N),N)=\sum_{i \geq 0} (-1)^{i} (i+1)\lambda^{i}Z_{i+1} (\bm{X},N),
\end{equation}
which turns out to be unbiased by construction. See Appendix \ref{sec:comparingE} for the determination of \eqref{eq:comparingE}. The estimator  $\hat{\tau}_{1}(\bm{X}(N),N)$ admits a natural interpretation as a nonparametric empirical Bayes estimator in the sense of \citet{Rob(56)}. More precisely, $\hat{\tau}_{1}(\bm{X}(N),N)$ is the posterior expectation of $\E[\tau_{1}(\bm{X},N,M)]$ with respect to an unknown prior distribution on the $p_{i}$'s that is estimated from the $Y_{j}(\mathbf{X},N)$. See Appendix \ref{sec:empirical_Bayes} for details. The next theorem legitimates the use of $\hat{\tau}_{1}(\bm{X}(N),N)$ as an estimator of $\tau_{1}(\bm{X},N,M)$, for $\lambda<1$, i.e. when the size of the unobserved population is less or equal than $n$, the size of the observed sample.

\begin{thm}\label{thm:var_lambda<1}
For any positive real numbers $x$ and $y$ let $\lfloor x \rfloor$ denote the integer part of $x$ and let $x\vee y$ denote the maximum between $x$ and $y$. If $\lambda<1$, for any $P \in \Pcr$, we get
\begin{equation}\label{eq:unbiased_T}
\E [\hat{\tau}_{1}(\bm{X}(N),N)]=\E [\tau_{1}(\bm{X},N,M)]=\sum_{j \geq 1}  np_j e^{-(\lambda +1)np_j}
\end{equation}
and 
\begin{equation}\label{bound_var_T}
\begin{split}
&\Var [\tau_{1}(\bm{X},N,M)-\hat{\tau}_{1}(\bm{X}(N),N)]\\
& \qquad\qquad\leq \Psi^2 (\lambda ) \E [Z_{\bar{1}} (\bm{X},N)]-\frac{\E [Z_{1} (\bm{X},N+M)]}{\lambda+1},
\end{split}
\end{equation}
where in \eqref{bound_var_T} we defined $\Psi (\lambda)= (j^{\ast} +1) \lambda^{j^{\ast}}$ such that $j^{\ast}=\lfloor (2 \lambda -1)/(1-\lambda) \rfloor \vee 0$.
\end{thm}

See Appendix \ref{sec:var_lambda<1} for the proof of Theorem \ref{thm:var_lambda<1}. According to Theorem \ref{thm:var_lambda<1}, for $\lambda<1$ one has $\E [\hat{\tau}_{1}(\bm{X}(N),N)]=\E [\tau_{1}(\bm{X},N,M)]$ and $\Var[\tau_{1}(\bm{X},N,M)-\hat{\tau}_{1}(\bm{X}(N),N)] \lesssim n$ upon noticing that $\E[Z_{\bar{1}} (\bm{X},N)]\leq \E [N]=n$. That is, in expectation, $\hat{\tau}_{1}(\bm{X}(N),N)$ approximate $\tau_{1}(\bm{X},N,M)$ to within $n$. Hence we formalize our considerations in the following.
\begin{cor}
Assume that $\lambda<1$, then  the nonparametric estimator $\hat{\tau}_{1}(\bm{X}(N),N)$ defined in \eqref{eq:estimator<1} satisfies
\begin{equation}
\Ecr_{\lambda , n} (\hat{\tau}_{1}(\bm{X}(N),N)) \lesssim \frac{1}{n}
\end{equation}
for any $n \geq 1$.
\end{cor}

This legitimates the use of $\hat{\tau}_{1}(\bm{X}(N),N)$ as an estimator of $\tau_{1}(\bm{X},N,M)$ under the hypothesis $\lambda<1$, which  unfortunately is a quite restrictive assumption within the framework of disclosure risk: indeed the size of the unobserved sample is usually much bigger than the size of the available one.
However the derivation of a variance bound for $\hat{\tau}_{1}(\bm{X}(N),N)$ is a crucial step for our study. Indeed it reveals that the assumption $\lambda<1$ is necessary to obtain a finite estimate of the variance. This variance issue of $\hat{\tau}_{1}(\bm{X}(N),N)$ is determined by the geometrically increasing magnitude of the coefficients $(i+1)(-\lambda)^{i}$. Indeed, as $\lambda\geq1$, the estimator $\hat{\tau}_{1}(\bm{X}(N),N)$ grows superlinearly as $(i+1)(-\lambda)^{i}$ for the largest $i$ such that $Z_{i+1} (\bm{X},N)>0$, thus eventually far exceeding $\tau_{1}(\bm{X},N,M)$ that grows at most linearly.
 This is the main reason why $\hat{\tau}_{1}(\bm{X}(N),N)$ become useless for $\lambda\geq1$, thus requiring an adjustment via suitable smoothing techniques. 
Hereafter we follow ideas originally developed by \citet{Goo(56)}, \citet{Efr(76)} and \citet{Orl(17)} for nonparametric estimators of the number of unseen species. Specifically, we propose a smoothed version of $\hat{\tau}_{1}(\bm{X}(N),N)$ by truncating the series \eqref{eq:estimator<1} at an independent random location $L$ and 
averaging over the distribution of $L$, i.e.,
\begin{align}\label{eq:estimator>1}
\hat{\tau}_{1}^L(\bm{X}(N),N)&= \E_L \left[ \sum_{i =1}^L (-1)^{i} (i+1) \lambda^i Z_{i+1} (\bm{X},N)\right]\\
&\notag=\sum_{i \geq 0} (-1)^{i} (i+1) \lambda^i \P (L \geq i) Z_{i+1} (\bm{X},N).
\end{align}
For any $\lambda\geq1$, as the the index $i$ in \eqref{eq:estimator>1} increases, the tail probability $\P[L\geq j]$ compensate for the exponential growth of $(i+1)(-\lambda)^{i}$, thereby stabilizing the variance. In the next theorem we show that for $\lambda\geq1$ the estimator $\hat{\tau}_{1}^L(\bm{X}(N),N)$ is biased for $\E[\tau_{1}(\bm{X},N,M)]$, and we provide a bound for the MSE of $\hat{\tau}_{1}(\bm{X}(N),N)$.

\begin{thm} \label{thm:var+bias_L}
Suppose that $\lambda \geq 1$, then $\hat{\tau}_{1}^L(\bm{X}(N),N)$ is a biased estimator of $\E[\tau_{1}(\bm{X},N,M)]$ with
\begin{equation}\label{eq:bias>1}
\begin{split}
&\E [\hat{\tau}_{1}^L(\bm{X}(N),N)] \\
& \qquad= \E [\tau_{1}(\bm{X},N,M)]+\sum_{j\geq 1} e^{-p_j n(\lambda +1)} p_j n \int_{0}^{\lambda np_j} e^s \E_{L} \left[ \frac{(-s)^L}{L!} \right] \D s.
\end{split}
\end{equation}
and
\begin{equation}
\label{eq:MSE_bound}
\begin{split}
&\MSE [\hat{\tau}^L_{1}(\bm{X}(N),N)] \\
& \qquad\leq \left( \sum_{j\geq 1} e^{-p_j n(\lambda +1)} p_j n \int_{0}^{\lambda np_j} e^s \E_{L} 
\left[ \frac{(-s)^L}{L!} \right] \D s \right)^2 \\
& \qquad\qquad + (\E_L [(L+1) \lambda^L])^{2} \E [Z_{\bar{1}} (\bm{X},N)] -\frac{\E [Z_{1} (\bm{X},N+M)]}{\lambda+1}  .
\end{split}
\end{equation}
\end{thm}

Choosing different smoothing distributions for the random variable $L$ yields different estimators for $\tau_{1}(\bm{X},N,M)$.
Following \citet{Orl(17)}, three possible choices for the distribution of $L$ are the following: i) a Poisson distribution with parameter $\beta >0$; ii) a Binomial distribution with parameter $(x_0, (\lambda+1)^{-1})$; iii) a Binomial distribution with parameter $(x_{0},2/(\lambda+2))$. In particular, it can be shown that the choice of the Binomial distribution with parameter $(x_0, (\lambda+1)^{-1})$ corresponds to the truncation at the point $x_{0}$ of the Euler transformation of the estimator \eqref{eq:estimator<1}. To choose the parameter $\beta$ of the Poisson distribution and the parameter $x_{0}$ of the Binomial distribution, one should look for $\tilde{\beta}$ and $\tilde{x}_{0}$ which minimizes the MSE bound \eqref{eq:MSE_bound}.
Once the values of $\tilde{\beta}$ and $\tilde{x}_{0}$ are determined explicitly, we are able to obtain limit of predictability for $\hat{\tau}_{1}^{L}(\bm{X}(N),N)$. That is,  for some $\delta >0$ we are able to specify the maximum value of the sampling fraction $\lambda$ for which $\Ecr_{\lambda , n} (\hat{\tau}_{1}^{L}(\bm{X}(N),N))< \delta$. This gives a provable (performance) guarantee for the estimation of $\tau_{1}(\bm{X},N,M)$ in terms of the sampling fraction $\lambda$. 
The next proposition specifies the limit of predictability for the estimator under the choice of   a Poisson distribution with parameter $\beta$ for the smoothing distribution $L$. 

\begin{prp} \label{prp:poisson}
Let $L$ be a Poisson random variable with parameter $\beta$. Then
\begin{equation} \label{eq:MSE_poisson}
\MSE [\hat{\tau}^L_{1}(\bm{X}(N),N)] \leq e^{-2\beta}n^2 + ne^{2\beta (2\lambda-1)}
\end{equation}
whose upper bound is minimized when
\begin{displaymath}
\tilde{\beta}=\frac{1}{4\lambda}\log \left( \frac{n}{2\lambda-1}\right).
\end{displaymath}
for any $\lambda\geq1$. Moreover, if $L$ is a Poisson random variable with parameter $\tilde{\beta}$ then
\begin{equation}
\label{eq:NMSE_poisson}
\Ecr_{n, \lambda} (\hat{\tau}^L_{1}(\bm{X}(N),N))  \leq \frac{A(\lambda)}{n^{1/(2\lambda)}},
\end{equation}
and for any $\delta\in(0,1)$
\begin{equation}
\label{eq:predictability_poisson}
\lim_{n \to +\infty} \frac{\max \left\{ \lambda : \; \Ecr_{n, \lambda} (\hat{\tau}^L_{1}(\bm{X}(N),N))\leq \delta \right\}}{\log (n)} \geq \frac{1}{2\log (A/\delta)}
\end{equation}
where $A(\lambda)$ is continuous in $[1, +\infty)$ with $\lim_{\lambda \to +\infty} A(\lambda)=1$ and $A= \max_{\lambda \geq 1} A(\lambda)< +\infty$.
\end{prp}
See Appendix \ref{sec:MSE_poiss} for the proof of Proposition \ref{prp:poisson}. Similar results hold true when $L$ is assumed to be a Binomial random variable: 
the derivation of these results follows  along similar lines as the proof of Proposition \ref{prp:poisson}. Hence we state the following result in presence of a Binomial smoothing without proof.

\begin{prp} \label{prp:bino2}
Let $L$ be a Binomial random variable with parameter $(x_{0},2/(\lambda+2))$. Then
\begin{equation}\label{eq:MSE_bino2}
\begin{split}
&\MSE [\hat{\tau}^L_{1}(\bm{X}(N),N)] \\
& \qquad\leq n 
\left( \frac{\lambda}{\lambda+2} \right)^{2 x_0} \left[3^{10 x_0/3}  +n \left(  \frac{\lambda}{2(\lambda+1)}\right)^2\right]
\end{split}
\end{equation}
whose upper bound is minimized when
\begin{displaymath}
\tilde{x}_0 = \left\lfloor \frac{3}{10} \log_3 \left( n \frac{\lambda^2 }{(\lambda+1) (\lambda^2 (3^{10/3}-1)-4 \lambda -4)} \right)   \right\rfloor
\end{displaymath}
for any $\lambda\geq1$. Moreover, if $L$ is a Binomial random variable with parameter $(\tilde{x}_{0},2/(\lambda+2))$ then
\begin{equation}\label{eq:NMSE_bino2}
\Ecr_{n, \lambda} (\hat{\tau}^L_{1}(\bm{X}(N),N))  \leq  \frac{C(\lambda)}{n^{3\log_3 (1+2/\lambda)/5}},
\end{equation}
and for any $\delta\in(0,1)$
\begin{equation}
\label{eq:predictability_binomial2}
\lim_{n \to +\infty} \frac{\max \left\{ \lambda : \;  \Ecr_{n, \lambda} (\hat{\tau}^L_{1}) \leq \delta\right\}}{\log(n)} \geq  
\frac{6}{5\log (3) \log (C/\delta)}
\end{equation}
where $C(\lambda)$ is continuous in $[1, +\infty)$ with $\lim_{\lambda \to +\infty} C(\lambda)=1$ and $C= \max_{\lambda \geq 1} C(\lambda)$.
\end{prp}

\section{Optimality of the proposed estimators}

In Section \ref{sec:theory} we have defined two different estimators of $\tau_1(\bm{X},N,M)$ providing guarantees of their performance as $n \to +\infty$ in terms of the NMSE. We have already remarked that the case $\lambda \geq 1$ is the most interesting one for estimating the disclosure risk $\tau_1(\bm{X},N,M)$, indeed the fraction of the unobserved sample $\lambda$ is usually much larger than $1$. Throughout the section we assume that $\lambda\geq 1$
and we prove that the proposed estimator $\hat{\tau}_1^L(\bm{X}(N),N)$ is essentially
optimal.
More precisely  we determine a lower bound for the best worst--case NMSE, defined by
\begin{equation} \label{eq:minimax}
\Ecr(\lambda , n) := \inf_{\hat{\rho}_1}\Ecr_{\lambda , n} (\hat{\rho}_1(\bm{X}(N),N))
\end{equation}
where the infimum in the previous definition runs over all possible estimators
of $\tau_1(\bm{X},N,M)$. We will then see that the determined
  lower bound essentially matches with the upper bound \eqref{eq:NMSE_poisson}.
In the sequel we  refer to 
$\Ecr(\lambda , n) $ as the \textit{(normalized) minimax risk}. 

The theorem we are going to state below provides us with a lower bound for the minimax risk.
\begin{thm}\label{thm:1}
  Assume that $1+\lambda > e^2$. Then, there exists a universal constant $K > 0$  such that for any $n$ sufficiently big we have
  \begin{equation} \label{eq:lower_geenral}
   \Ecr(\lambda , n)  
    \geq
K \cdot \left\{ \begin{array}{ll}
 1 & \text{if } \lambda +1 > \log (n)\\
 \frac{1+\lambda}{\log (n)} \left( \frac{\sqrt{\log(n)}}{n (1+\lambda)} \right)^{e^2/(1+\lambda)}
 & \text{if } \lambda + 1\leq \log (n)
 \end{array}
  \right. 
  \end{equation}
\end{thm}
From Theorem \ref{thm:1}, it is clear that the minimax risk goes
to zero if $\lambda+1 = o( \log (n))$ and the rate is provided by the following Corollary.
\begin{cor} \label{cor:lower_bound}
  Assume that $1+\lambda > e^2$, then there exist universal constants $c >0 $ and $c'>0$ such that for any $n$ sufficiently large 
  \begin{equation}
  \label{eq:lower_bound}
     \Ecr(\lambda , n)  
    \geq c \frac{1}{n^{c'/\lambda}}.
\end{equation}   
\end{cor}
Corollary \ref{cor:lower_bound} is an easy consequence of Theorem \ref{thm:1}, indeed, when $\lambda+1 >\log (n)$
the two lower bounds in \eqref{eq:lower_geenral}--\eqref{eq:lower_bound} are constants,  whereas if $\lambda +1 \leq \log (n)$ it is easy to observe that the leading term  in \eqref{eq:lower_geenral} (as $n \to +\infty$) is of order $1/n^{c'/\lambda}$ as in  \eqref{eq:lower_bound} for some $c'>0$.
Corollary \ref{cor:lower_bound} provides us with a lower bound for the NMSE of any estimator of the disclosure risk $\tau_1(\bm{X},N,M)$. The lower bound \eqref{eq:lower_bound} has an important implication:  without imposing any parametric assumption on the model, one can estimate  $\tau_1(\bm{X},N,M)$ with  vanishing NMSE all the way up to $\lambda \propto \log n$.
 It is then impossible to determine an estimator having provable guarantees (in
 terms of vanishing NMSE) when $\lambda = \lambda (n)$ goes to $+\infty$ much faster than $\log (n)$,
 as a function of $n$. By the limit of predictability \eqref{eq:predictability_poisson} determined for 
the estimator $\hat{\tau}_1^L(\bm{X}(N),N)$, we conclude that the proposed estimator  is optimal, because its limit of predictability matches (asymptotically) with its maximum possible value $\lambda \propto \log (n)$.

\subsection{Guideline for the proof of Theorem \ref{thm:1}}
\label{sec:proof-lower-bound}

In the present section we provide the main ingredients for the proof of Theorem \ref{thm:1}, technical results and related proofs are deferred to the Appendix.
In the sequel we will write $\E_P^n$ to make explicit the dependence of the expected value w.r.t. $P$ and the parameter $n$ of the Poisson random variable $N$.

The starting point for the proof of Theorem \ref{thm:1} is the next Lemma \ref{lem:3}, which is an interesting result in its own right and will help a lot in the proof of Theorem \ref{thm:1}. Remark that the definition of the minimax risk in  \eqref{eq:minimax}
allows for estimators depending on the whole sample $\bm{X}(N)$, while
$\tau_1(\bm{X},N,M)$ depends only on the frequencies $\bm{Y}(\bm{X},N+M)$ and
$\bm{Y}(\bm{X},N)$. Thus, we feel like there should be no gain of
information in using estimators depending on $\bm{X}(N)$ over estimators
depending only on the frequencies $\bm{Y}(\bm{X}, N)$. This is made
formal in the next lemma, proved in Section \ref{sec:proof-lem:3}. Note that this is
convenient since $(\bm{X},k) \mapsto \bm{Y}(\bm{X},k)$ is nicely
distributed under the Poisson model.
\begin{lem}
  \label{lem:3}
The following equality is true
  \begin{equation*}
    \Ecr(\lambda , n) %
    = \inf_{\hat{\rho}}\sup_{P\in
      \Pcr} n^{-2}\E_P^n[(\tau_1(\bm{X},N,M) -
    \hat{\rho}(\bm{Y}(\bm{X},N))^2],
  \end{equation*}
where the infinimum in the previous equation is understood to be taken with respect
  to all measurable maps $\hat{\rho}: \N^{\N} \rightarrow\R$.
\end{lem}

\par The next step is to use Jensen's inequality to deduce that
\begin{align*}
  \Ecr(\lambda , n) %
  & = \inf_{\hat{\rho}}\sup_{P\in \Pcr}n^{-2} \E_P^n[ \E_P^n[(
    \tau_1(\bm{X},N,M) - \hat{\rho}(\bm{Y}(\bm{X},N)))^2  \mid
    \bm{Y}(\bm{X},N)] ]\\
  &\geq \inf_{\hat{\rho}}\sup_{P\in \Pcr} n^{-2} \E_P^n[ (
    \E_P^n[\tau_1(\bm{X},N,M) \mid \bm{Y}(\bm{X},N)] -
    \hat{\rho}(\bm{Y}(\bm{X},N)))^2  ]
\end{align*}
Note that there is no explicit dependency on $\bm{X}$ and $M$ anymore in
the last display, but only on the random variable
$(\bm{X},N) \mapsto \bm{Y}(\bm{X},N)$ which, under $P$, is
distributed as an infinite vector of independent Poisson random variables with
parameters $(n p_1,n p_2,\dots)$. Besides observe also that
$N = \sum_{j \geq 1}Y_j(\bm{X},N)$. For the sake of notational simplicity, in the sequel 
$\bm{Y}$ will stand for the random variable
$(\bm{X},N) \mapsto \bm{Y}(\bm{X},N)$, and we also let
\begin{align*}
  \tilde{\tau}_1(\bm{Y},P,n)%
  &:= \E_P^n[\tau_1(\bm{X},N,M) \mid \bm{Y}(\bm{X},N)]\\
  &= \sum_{j\geq 1}\ind{\{Y_j(\bm{X},N) =1\}}%
    \E_P^n[\ind{\{Y_j(\bm{X},N+M ) - Y_j(\bm{X},N) = 0\}} \mid
    \bm{Y}(\bm{X},N) ].
\end{align*}
Remark that $(Y_j(\bm{X},N+M) - Y_j(\bm{X},N):\, j\in \N)$ is
independent of $\bm{Y}(\bm{X},N)$ and is a collection of independent
Poisson random variables with intensities
$(\lambda n p_j :\, j\in \N)$. Henceforth, we get
\begin{equation}
  \label{eq:7}
  \tilde{\tau}_1(\bm{Y},P,n)%
  = \sum_{j\geq 1}e^{-\lambda n p_j}\ind{\{Y_j (\bm{X},N) = 1\}},
\end{equation}
and besides
\begin{equation}
  \label{eq:10}
  \Ecr(\lambda , n) %
  \geq \inf_{\hat{\rho}}\sup_{P\in \Pcr} n^{-2}\E_P^n[(
  \tilde{\tau}_1(\bm{Y},P,n) - \hat{\rho}(\bm{Y}))^2 ].
\end{equation}

We now trade $\tilde{\tau}_1(\bm{Y},P,n)$ for its expectation. Let us introduce
$\bar{\tau}_1(P,n) := \E_P^n[ \tilde{\tau}_1(\bm{Y},P,n) ]$. Recall
that under $P$ the vector $\bm{Y}$ is distributed as independent Poisson
with parameters $(n p_1, np_2,\dots)$. Hence,
\begin{equation*}
  \bar{\tau}_1(P,n)%
  = \sum_{j\geq 1} e^{-\lambda n p_j}\E_P^n[\ind{\{Y_j (\bm{X},N) = 1\}}]%
  = n \sum_{j\geq 1}p_j e^{-(1+\lambda) n p_j}.
\end{equation*}
Similarly, for any $P \in \Pcr$,
\begin{equation}
  \label{eq:6}
  \Var(\tilde{\tau}_1(\bm{Y},P,n))%
  = \sum_{j\geq 1} n p_j e^{-(1+2\lambda) n p_j}\big\{1 - np_je^{- n p_j}
  \big\}%
  \leq n.
\end{equation}
Thus from \eqref{eq:10}, Young's inequality, we find that
\begin{align}
  \label{eq:11}
  \Ecr(\lambda , n) %
  &\geq \frac{1}{2} \inf_{\hat{\rho}}\sup_{P\in \Pcr}
    n^{-2}\E_P^n[(\bar{\tau}_1(P,n) - \hat{\rho}(\bm{Y}))^2 ] - n^{-1}.
\end{align}

The remainder of the proof mostly follows the reduction scheme used in
\citet{wu2016minimax,wu2015chebyshev} which consists on reducing the problem to
finding the best polynomial approximation (in uniform norm) to a suitable
function.

The first step of the reduction scheme is to trade $\Pcr$ in \eqref{eq:minimax} for a
slightly more convenient set. We let $S \in \N$ be the only integer
satisfying
\begin{equation}
  \label{eq:20}
  n(1+\lambda) \leq S \leq n(1+\lambda) + 1,
\end{equation}
and we also let, for some constant $c_0 > 0$ to be determined later,%
\begin{equation}
  \label{eq:21}
  \xi := (2c_0/e)\min\{(1+\lambda)\log n,\, \log^2 n\}.
\end{equation}
Then, for another constant $c_1 > 0$ and for $\varepsilon > 0$ to be determined
later, we define 
\begin{equation}
  \label{eq:22}
  \Pcr'
  :=
  \left\{\textstyle\sum_{k=1}^Sp_k
    \delta_k : \;  p_k \in [0,\xi S^{-1}],\ |\sum_{k=1}^S p_k - 1| \leq c_1
    \varepsilon / \xi \right\}.
\end{equation}
Remark that $\Pcr'$ contains measures that are not probability measures, and
hence it is not clear a priori that we can lower bound the supremum over
$\Pcr$ by the supremum over $\Pcr'$. The next proposition shows that it is
fine as long as $c_1 \varepsilon$ is not too large. Here and after, under
$P \in \Pcr'$, $\bm{Y}$ is understood as a vector of independent Poisson
random variables with intensities $(np_1,\dots, n p_S,0,\dots)$, with
$\sum_{j=1}^S p_j$ not necessarily equal to one, and $P \mapsto \bar{\tau}_1(P)$
is extended trivially from $\Pcr$ to $\Pcr'$ by letting
$\bar{\tau}_1(P,n) := n\sum_{j=1}^Sp_je^{-n(1+\lambda)p_j}$,
$P\in \Pcr'$. The next proposition is proved in Section \ref{sec:proof-pro:6}.

\begin{prp}
  \label{pro:6}
  Assume that $c_1\varepsilon = o(\xi)$ as $n\rightarrow \infty$ and let define
  $n' := n(1 + c_1\varepsilon/\xi)$. Then as $n \rightarrow \infty$,
  \begin{align*}
    \Ecr(\lambda , n) %
    &\geq  \frac{1}{4}\inf_{\hat{\rho}}\sup_{P \in \Pcr'}n^{-2}\E_P^{n'}[
      (\bar{\tau}_1(P,n) - \hat{\rho}(\bm{Y}))^2]%
      - (n^{-1} + 9 c_1^2 \varepsilon^2/2)\\
    &\geq \frac{\varepsilon^2}{4}\Big\{\inf_{\hat{\rho}}\sup_{P \in
      \Pcr'}\P_P^{n'}\big( |\bar{\tau}_1(P,n) - \hat{\rho}(\bm{Y})| > n
      \varepsilon \big) - 18c_1^2 \Big\} - n^{-1}.
  \end{align*}
\end{prp}

\par We are now in position to lower bound the risk by the Bayes risk. To do so,
we follow the prior construction of \citet{wu2016minimax,wu2015chebyshev}. For
some $L \in \N$ to be determined later, but satisfying $L \leq c_2\xi$ for
some constant $c_2 > 0$, we let $U$ and $V$ be two random variables taking
values in $[0,\xi S^{-1}]$ such that when $n$ is large enough,
\begin{gather*}
  \E[U^k] = \E[V^k]\quad \forall k \in \{0,\dots,L+1\},\\
  \E[U] = \E[V] = S^{-1},\quad%
  \Var(U) \leq \xi S^{-2},\quad \Var(V) \leq \xi S^{-2},\\
  \E[U e^{-n(1+\lambda)U}] \geq \E[V e^{-n(1+\lambda)V}]%
  + S^{-1} K\min\big\{1,\, \sqrt{\xi/L^2}\exp(- L^2/\xi) \big\}
\end{gather*}
The existence of such random variables is guaranteed by Theorem \ref{thm:4} for a
universal constant $K > 0$. Then we let $\bm{U} := (U_1,\dots,U_S)$,
respectively $\bm{V} := (V_1,\dots,V_S)$, be an independent vector of
i.i.d. copies of $U$, respectively $V$. Denoted by $\mathcal{M}(\N)$ the space of all measures on $\N$, we construct the following random variable
$Q : [0,1]^S \rightarrow \mathcal{M}(\N)$
such that $Q(\bm{U}) := \sum_{k=1}^S U_k \delta_k$. 
Then, from the Proposition \ref{pro:6} and Hölder's inequality, we find that
$\Ecr(\lambda,n)$ is bounded from below by $-n^{-1}$ plus
\begin{multline*}
  \frac{\varepsilon^2}{4}\Big\{%
  \inf_{\hat{\rho}}\Big(
  \frac{1}{2}\E\Big[\P_{Q(\bm{U})}^{n'}\big(
  |\bar{\tau}_1(Q(\bm{U}),n) - \hat{\rho}(\bm{Y})| > n\varepsilon
  \big) \mathbbm{1}_{\Pcr'}(Q(\bm{U})) \Big]\\
  + \frac{1}{2}\E\Big[ \P_{Q(\bm{V})}^{n'}\big( |
  \bar{\tau}_1(Q(\bm{V}),n) - \hat{\rho}(\bm{Y})| > n\varepsilon
  \big)\mathbbm{1}_{\Pcr'}(Q(\bm{V})) \Big]\Big) - 18c_1^2\Big\},%
\end{multline*}
which is in turn lower bounded by
\begin{multline*}
  \frac{\varepsilon^2}{4}\Big\{%
  \inf_{\hat{\rho}}\Big(%
  \frac{1}{2}
  \E\Big[\P_{Q(\bm{U})}^{n'}\big(
  |\bar{\tau}_1(Q(\bm{U}),n) - \hat{\rho}(\bm{Y})| > n\varepsilon
  \big) \Big]\\
  + \frac{1}{2}\E\Big[ \P_{Q(\bm{V})}^{n'}\big( |
  \bar{\tau}_1(Q(\bm{V}),n) - \hat{\rho}(\bm{Y})| > n\varepsilon
  \big) \Big] \Big)\\%
  - 18c_1^2%
  - \frac{1}{2}\P( Q(\bm{U}) \notin \Pcr')%
  - \frac{1}{2}\P( Q(\bm{V}) \notin \Pcr')%
  \Big\} - n^{-1}.%
\end{multline*}
The last display follows because we don't have
$Q(\bm{U})$ nor $Q(\bm{V})$ almost-surely in $\Pcr'$, but it is clear
that the strong law of large numbers implies they should be concentrated on
$\Pcr'$. Formally, an application of Bernstein's inequality (see
Section \ref{sec:proof-pro:2} below) leads to the following proposition.
\begin{prp}
  \label{pro:2}
  Assume that $c_1\varepsilon = o(\xi)$ as $n\rightarrow \infty$. Then, there
  exists a constant $C > 0$, depending only on $c_1$, such that for $n$ large
  enough,
  \begin{equation*}
    \varepsilon^2%
    \geq \frac{C \xi^3}{n(1+\lambda)}
    \Longrightarrow%
    \P\big( Q(\bm{U}) \notin \Pcr' \big)%
    \leq c_1^2.
  \end{equation*}
\end{prp}

Thus under the conditions of Proposition \ref{pro:2}, we get
\begin{multline}
  \label{eq:3}
  \Ecr(\lambda , n) %
  \geq \frac{\varepsilon^2}{4} \Big\{%
  \inf_{\hat{\rho}}\Big( \frac{1}{2}\E\Big[\P_{Q(\bm{U})}^{n'}\big(
  |\bar{\tau}_1(Q(\bm{U}),n) - \hat{\rho}(\bm{Y})| > n\varepsilon
  \big) \Big]\\
  + \frac{1}{2}\E\Big[ \P_{Q(\bm{V})}^{n'}\big( |
  \bar{\tau}_1(Q(\bm{V}),n) - \hat{\rho}(\bm{Y})| > n\varepsilon
  \big) \Big] \Big) - 19 c_1^2\Big\} - n^{-1}.
\end{multline}
We now wish to trade $\bar{\tau}_1(Q(\bm{U}),n)$ and
$\bar{\tau}_1(Q(\bm{V}),n)$ for their expectations in the last
equation. Intuitively, this should not be problematic since they are sums of i.i.d.
random variables, they should concentrate near their expectations for $S$ large
enough. We made this formal using a Hoeffding argument in the next proposition,
proved in Section \ref{sec:proof-pro:7}.

\begin{prp}
  \label{pro:7}
  Let everything as above. Then,
  \begin{equation*}
    \varepsilon^2%
    \geq \frac{2\xi \log(1/c_1^2)}{n(1+\lambda)}
    \Longrightarrow%
    \P\big( | \bar{\tau}_1(Q(\bm{U}),n) -
    \E[\bar{\tau}_1(Q(\bm{U}),n)]| >  n\varepsilon /2 \big) \leq 2c_1^2.
  \end{equation*}
\end{prp}

Obviously, Proposition \ref{pro:7} is also true for $\bar{\tau}_1(Q(\bm{V}),n)$. We
now assume that the conditions of Propositions \ref{pro:6}, \ref{pro:2} and \ref{pro:7} are met, and thus we
obtain from \eqref{eq:3} that
\begin{multline*}
  \Ecr(\lambda , n) %
  \geq \frac{\varepsilon^2}{4} \Big\{%
  \inf_{\hat{\rho}}\Big( \frac{1}{2}\E\Big[\P_{Q(\bm{U})}^{n'}\big(
  |\E[\bar{\tau}_1(Q(\bm{U}),n)] - \hat{\rho}(\bm{Y})| > n\varepsilon/2
  \big) \Big]\\
  + \frac{1}{2}\E\Big[ \P_{Q(\bm{V})}^{n'}\big( |
  \E[\bar{\tau}_1(Q(\bm{V}),n)] - \hat{\rho}(\bm{Y})| > n\varepsilon/2
  \big) \Big] \Big) - 21 c_1^2\Big\} - n^{-1}.
\end{multline*}
Now remark that,
\begin{equation*}
  \E[\bar{\tau}_1(Q(\bm{U}),n)]%
  = n S \E[U \exp\{-n(1+\lambda)U\}],
\end{equation*}
and besides observe that whenever $n$ is large enough, we will have
$C\xi^2 \geq 2\log(1/c_1^2)$. We furthermore assume that $\varepsilon^2$
satisfies
\begin{equation}
  \label{eq:8}
  \max\Big\{\frac{8}{n c_1^2},\,\frac{C\xi^3}{n(1+\lambda)} \Big\}%
  \leq \varepsilon^2 \leq%
  K^2 \min\big\{1,\, (\xi/L^2)\exp(-2L^2/\xi) \big\}.
\end{equation}
It is not clear yet that \eqref{eq:8} can be satisfied, we claim this is the case
and we delay the proof of the claim at the end of the section. When the claim is
true, we pick $\varepsilon$ to be equal to the r.h.s. of \eqref{eq:8}. Then, the
previous computations and the classical Le Cam method with two hypothesis imply
that
\begin{align*}
 & \Ecr(\lambda , n) \\
  &\quad\geq \frac{\varepsilon^2}{8} \Big\{%
  1 - \TV\Big(\otimes_{j=1}^S \E[{\rm Poiss}(n' U_j)], \otimes_{j=1}^S
  \E[{\rm Poiss}(n' V_j)] \Big) - 42 c_1^2 \Big\} - \frac{1}{n}\\
  &\quad\geq%
    \frac{\varepsilon^2}{8} \Big\{%
  1 - S \cdot \TV\Big(\E[{\rm Poiss}(n' U)], \E[{\rm Poiss}(n' V)]
    \Big) - 43c_1^2 \Big\}.
\end{align*}

\par We now make explicit our choice for the value of $L$. We need to choose it
small enough such that $\E[U \exp\{-n(1+\lambda)U\}]$ and
$\E[V \exp\{-n(1+\lambda)V\}]$ are maximally separated, but also large enough
such that the previous display is non-negative. 
  For a constant $c_3 > 0$ to be chosen accordingly later, and for
  $A(\lambda,n) > 0$ solution of%
\begin{equation*}
  A(\lambda,n)\log A(\lambda,n)%
  = c_0^{-1} + c_0^{-1}\frac{\log(1+\lambda) - (1/2)\log\log(n) +
    \log(c_3)}{\log(n)},
\end{equation*}
we pick $L$ to be the smallest integer
satisfying the bound,
\begin{equation}
  \label{eq:23}
  L \geq%
  \begin{cases}
    2 c_0 \log(n) &\mathrm{if}\ 1+\lambda > \log(n),\\
    c_0A(\lambda,n) \log(n) &\mathrm{if}\ 1+\lambda \leq \log(n).
  \end{cases}
\end{equation}
Remark that this choice ensure that $L \leq c_2 \xi$ for some constant
$c_2 > 0$, as requested. Without loss of generality we further assume that
$L+2 \leq c_2\xi$. We are then able to state the following propositon, which will be proved in
Section \ref{sec:proof-pro:8}.
\begin{prp}
  \label{pro:8}
  Assume that $c_1c_2 \varepsilon \leq 1$. Then, the constant $c_3 > 0$ can be
  chosen depending only on the choice of $c_0$ and such that for $n$ large
  enough,
  \begin{equation*}
    S \TV\Big(\E[{\rm Poiss}(n' U)], \E[{\rm Poiss}(n' V)]\Big) \leq \frac{1}{2}.
  \end{equation*}
\end{prp}
Now observe that $\varepsilon \leq K$  by \eqref{eq:8}, thus choosing
$
c_1 = \min\{1/\sqrt{172},\, 1/(c_2K)\},
$
 we get
$c_1c_2\varepsilon \leq 1$ and $43c_1^2 \leq 1/4$. Together with the last
proposition, this implies that $\Ecr(\lambda , n)  \geq \varepsilon^2/32$, at least
if $\varepsilon$ can be chosen to satisfy Equation \eqref{eq:8}.

\par We now prove that \eqref{eq:8} is satisfied. When $1+\lambda > \log (n)$, then the r.h.s. of \eqref{eq:8} is always greater
than a constant, while the the l.h.s. is always smaller than $O(1/n)$. Thus in this
case we take $\varepsilon^2$ to be a suitable constant, and \eqref{eq:8} is
satisfied, giving $\Ecr(\lambda , n)  \gtrsim 1$.

\par When $1 + \lambda \leq \log (n)$, then there is always a constant
$K' > 0$ such that the r.h.s. of \eqref{eq:8} satisfies,
\begin{align*}
  &\min\big\{1,\, \sqrt{\xi/L^2}\exp(-L^2/\xi) \big\}
  \geq K' \frac{\sqrt{\xi}}{L}\exp(-L^2/\xi)\\
  &\qquad\qquad\geq \frac{K'\sqrt{2/e}}{\sqrt{c_0}A(\lambda,n)}%
    \frac{ \sqrt{1+\lambda}}{\sqrt{\log(n)}}%
    \exp\Big\{- \frac{e c_0A(\lambda,n)^2}{2(1+\lambda)} \log(n) \Big\}.
\end{align*}
We now make explicit our choice for $c_0$. We pick $c_0 = 1/e$, and it can be
seen that this choice is the one which asymptotically minimizes the product
$c_0 A(\lambda,n)^2$. Then, the Proposition \ref{pro:1} below, proven in
Section \ref{sec:proof-pro:1}, yields to the following bound, for a universal constant
$K'' > 0$, and for $n$ large enough,
\begin{equation*}
  \min\big\{1,\, \sqrt{\xi/L^2}\exp(-L^2/\xi) \big\}%
  \geq%
  \frac{K''\sqrt{1+\lambda}}{\sqrt{\log(n)}}%
  \Big( \frac{\sqrt{\log(n)}}{n(1+\lambda)} \Big)^{\frac{e^2/2}{1+\lambda}}.
\end{equation*}
The last equation then shows that whenever $1+\lambda > e^2$, the r.h.s. of
\eqref{eq:8} gets larger than the l.h.s. when $n$ is large enough. This concludes
the proof of the theorem.

\begin{prp}
  \label{pro:1}
  Let $c_0 = 1/e$. Then whenever $1+\lambda \leq \log(n)$ we have
  $A(\lambda,n) = e + o(1)$ as $n\rightarrow \infty$. Furthermore when
  $1+\lambda \leq \log(n)$, as $n \rightarrow \infty$,
  \begin{equation*}
    c_0A(\lambda,n)^2 \log(n)%
    \leq  e \log(n)%
    + e\log \frac{c_3(1+\lambda)}{\sqrt{\log(n)}}%
    +o(1).
  \end{equation*}
\end{prp}

\section{Numerical illustrations}

We present an illustration on synthetic data of the estimators introduced in
Section \ref{sec:theory}. We also consider other estimators of $\tau_{1}$ that
have been proposed in the literature of disclosure risk assessment: i) two
parametric empirical Bayes estimators of $\tau_{1}$ proposed by \citet{Bet(90)}
and \citet{Ski(94)}; ii) a naive nonparametric estimator of $\tau_{1}$; iii) a
Bayesian nonparametric estimator of $\tau_{1}$ proposed by \citet{Sam(98)}. A
common feature of these estimators, as well as our class of nonparametric
estimators, is that they rely on modeling the random partition induced by the
cross-classified sample records. More recent approaches, not considered here,
focus on modeling associations among identifying variables by log-linear models,
local smoothing polynomials and hierarchical latent models. In particular, the
Bayesian hierarchical semiparametric models of \citet{Car(15)} and
\citet{Car(18)} show a remarkable better performance than models
for random partitions, at the cost of an increasing computational effort for the need of Markov chain Monte Carlo methods for posterior approximation.

The approach of \citet{Bet(90)} is a parametric empirical Bayes approach in the sense of \citet{Efr(73)}. It relies on the following modeling assumption for the cells' frequencies of the population: $Y_j (\bm{X},\bar{n}) \sim {\rm Poiss} (\bar{n} p_j)$, where $\bar{n}$ is the size of the entire population. \citet{Bet(90)} also assumed a Gamma prior distribution over the probabilities associated to each cell, namely $p_j \sim {\rm Gam} (\alpha, \beta)$. One should specify the $p_j$'s under the condition $\sum_{j=1}^{K_{\bar{n}}} p_j =1$, however, for the sake of simplicity, \citet{Bet(90)} assumed that $\sum_{j =1}^{K_{\bar{n}}}\E [p_j] =1$, which is tantamount to saying that $\alpha= 1/(K_{\bar{n}}\beta)$. Under these modeling assumptions, \citet{Bet(90)} proposed an estimator of the expected value of total number  $T_1 (\bm{X},\bar{n})$ of population uniques, i.e., 
\begin{equation}
\label{eq:T1}
T_1 (\bm{X},\bar{n}) := \sum_{j =1}^{K_{\bar{n}}} \ind{\{ Y_j (\bm{X},\bar{n})=1\}}.
\end{equation}
Under the above Poisson-Gamma model, $\E [T_1 (\bm{X},\bar{n}) ]= \bar{n} (1+\bar{n} \beta)^{-(1+\alpha)}$, which depends on the parameters $\alpha$ and $\beta$, with the condition $\alpha=1/(K \beta)$. Parameters can be easily estimated via maximum likelihood, as we have done in the subsequent numerical experiments. If $K_{\bar{n}}$ is not available, \citet{Bet(90)} suggested to estimate $K_{\bar{n}}$ assuming  a uniform distribution over the cells, hence 
\[
\hat{K}_{\bar{n}} = \frac{\bar{n}K_n}{\sum_{j=1}^{K_{n}}   \ind{\{ Y_j (\bm{X},n)=1\}} },
\]
where $n$ is the size of the observed sample and $K_n$ stands for the number of distinct cells dictated by the sample of size $n$. If $\hat{\alpha}$ and $\hat{\beta}$ denote the maximum likelihood estimators of $\alpha$ and $\beta$, respectively, then an estimator of $T_1 (\bm{X},\bar{n})$ is $\hat{T}_1 = \bar{n} (1+\bar{n} \hat{\beta})^{-(1+\hat{\alpha})}$. \citet{Bet(90)} then suggested a corresponding estimator of $\tau_1$ as the sample portion of $\hat{T}_1$. More precisely, they proposed
\begin{equation}
\label{eq:Bet}
\hat{\tau}_1^{B} =\frac{n}{\bar{n}}\hat{T}_1 =n (1+\bar{n} \hat{\beta})^{-(1+\hat{\alpha})}.
\end{equation}
as an estimator of $\tau_{1}$.
\citet{Ski(94)} improved the estimator \eqref{eq:Bet}. In particular, still under the Poisson-Gamma model, they considered directly the problem of estimating $\tau_{1}$. In particular, they proposed the following estimator
\begin{equation}
\label{eq:Ski}
\hat{\tau}_1^S := K_n \left( \frac{1+\bar{n}\hat{\beta}}{1+n\hat{\beta}} \right)^{-(1+\hat{\alpha})},
\end{equation}
where the prior parameters $\alpha$ and $\beta$ can be estimated via maximum likelihood. The estimators proposed in Section \ref{sec:theory}, due to their nonparametric empirical Bayes interpretation in the sense of \citet{Rob(56)}, may be considered as the natural nonparametric counterparts of the empirical Bayes estimator \eqref{eq:Ski}.

Besides parametric estimators of $\tau_{1}$, we also consider two nonparametric estimators. A naive nonparametric estimator of $\tau_{1}$ relies on the intuition that a natural estimator of $\tau_{1}$ is the sampling fraction, with respect to the population, of the number of sample uniques. This estimator was first discussed in \citet{Bet(90)} and \citet{Ski(94)}, and it is defined as follows
\begin{equation}
\label{eq:naive}
\hat{\tau}^{\Ncr}_1 := Z_1(\bm{X},n)\frac{n}{\bar{n}}.
\end{equation}
\citet{Sam(98)} exploits Bayesian nonparametric ideas, and in particular a Dirichlet process prior (\citet{Fer(73)}) on the $p_{j}$'s to derive a smoothed version of the naive estimator \eqref{eq:naive}. In particular, \citet{Sam(98)} suggested the following estimator
\begin{equation}
\label{eq:samuels}
\hat{\tau}^{\Dcr}_1 := Z_1(\bm{X},n)\frac{n+\vartheta -1}{\bar{n}+\vartheta-1},
\end{equation}
where $\vartheta$ is the concentration parameter of the Dirichlet process prior. It is well-known (see, e.g. \citet{Fer(73)}) that the maximum likelihood estimator of $\vartheta$ can be obtained by solving, with respect to $\vartheta$, the equation $K_n = \sum_{1\leq j\leq n-1}\vartheta/(\vartheta+j)$.

We generate synthetic tables with $C$ cells, where $C=3 \cdot 10^5$ for Table \ref{tab:simulated_example_s}, $C=6 \cdot 10^5$ for Table \ref{tab:simulated_example_m} and $C=9 \cdot 10^5$ for Table \ref{tab:simulated_example_l}. For any choice of the size $C$, the true probabilities $(p_{j})_{j\geq1}$ of cells have been generated according to different types of distributions: the Zipf distribution, i.e., $p_{j}\propto j^{-s}$ for some $s>0$, the uniform distribution over the total number of cells and the uniform Dirichlet distribution. For all the simulated scenarios we have considered  a population of size $\bar{n}= 10^6$ and a sample of $n=10^{5}$ individuals from it. Each column of Tables \ref{tab:simulated_example_s}--\ref{tab:simulated_example_l} corresponds to a different choice of the distribution over the cells' probabilities. From the left to right: the Zipf distribution with parameter $s=0.2, 0.5,  0.8, 1$, the uniform distribution, the uniform Dirichlet distribution with parameter $\beta=0.5, 1$. In the first row of each table we have reported the true values of the disclosure index, while the other rows contain the estimates obtained with: i) the nonparametric estimator with Binomial smoothing $\hat{\tau}_1^{L_b}$, see Proposition \ref{prp:bino2}; ii) the nonparametric estimator with Poisson smoothing $\hat{\tau}_1^{L_p} $, see Proposition \ref{prp:poisson}; iii) the naive nonparametric estimator $\hat{\tau}^{\Ncr}_1$; iv) the Bayesian nonparametric estimator $\hat{\tau}^{\Dcr}_1$; v) the parametric empirical Bayes estimator $\hat{\tau}_1^{B}$; vi) the parametric empirical Bayes estimator $\hat{\tau}_1^{S}$
All experiments are averaged over 100 iterations. The best estimates in each simulated scenarios are displayed in bold. \\

\begin{table}[h!]
\begin{center}
\begin{tabular}{{|c||ccccccc|}}\hline
{}  & {Zipf $0.2$} & {Zipf $0.5$} & {Zipf $0.8$} & {Zipf $1$}  & {Uniform} & {Dirichlet $0.5$} & {Dirichlet $1$} \\ \hline\hline
True $\tau_1$   &  2868   &  4651  & 7537 & 7313 & 2579  & 4151  & 4608\\\hline
$\hat{\tau}_1^{L_b}$ & \textbf{5671} &9322  &  12790 & 10373  & \textbf{5582} & 2814 & 3118 \\[0.16cm]
$\hat{\tau}_1^{L_p}$ & 21086 & 20451 & 18221  &  12744 & 21656 & 8205 &  11855\\[0.16cm]
$\hat{\tau}_1^{\Ncr}$ & 6413 & \textbf{5613} & 3810   & 2157  & 6511  & \textbf{4232} & \textbf{5111}\\[0.16cm]
$\hat{\tau}_1^{\Dcr}$ & 18461 & 12471 & \textbf{5421}  & 2459 & 19303 &  7498  & 10741 \\[0.16cm]
$\hat{\tau}_1^{B}$ & 30554  &  27271 & 13848 &\textbf{ 5358}  & 30847   & 22820  &  26351 \\[0.16cm]
$\hat{\tau}_1^{S}$ &  28702 & 23621 & 9670  & 3043 &  29187  & 17665   & 22306\\\hline
\end{tabular}
\begin{flushleft}
\caption{Estimators of $\tau_1$ for several simulated scenarios, when the size of the table is $C=3\cdot 10^5$.}
\label{tab:simulated_example_s}
\end{flushleft}
\end{center}
\end{table}

\clearpage

\begin{table}[h!]
\begin{center}
\begin{tabular}{{|c||ccccccc|}}\hline
{}  & {Zipf $0.2$} & {Zipf $0.5$} & {Zipf $0.8$} & {Zipf $1$}  & {Uniform} & {Dirichlet $0.5$} & {Dirichlet $1$} \\ \hline \hline
True $\tau_1$   &   16020  & 16819  & 16095  &  11401   & 15947 & 9857 & 12451\\\hline
$\hat{\tau}_1^{L_b}$ & 32254  & 27585  & 21216  & \textbf{13194}  & 33426 &  \textbf{9190} & \textbf{14805}\\[0.16cm]
$\hat{\tau}_1^{L_p}$ & 49883  & 42310  & 29344 & 16847  & 51096  & 22932  & 32146\\[0.16cm]
$\hat{\tau}_1^{\Ncr}$ &  \textbf{7625}  & 6860  & 4642  & 2534 & \textbf{7693} & 5899  & 6670\\[0.16cm]
$\hat{\tau}_1^{\Dcr}$ & 32567 & \textbf{21009} & 7380  & 2968 & 33860  & 14577  &  20337\\[0.16cm]
$\hat{\tau}_1^{B}$ & 33594  & 31155    & \textbf{17152} & 6380  & 33763 & 28795  & 31114\\[0.16cm]
$\hat{\tau}_1^{S}$ & 34070   & 29703  & 13032  & 3776  & 34321 & 25900  & 29634\\\hline
\end{tabular}
\begin{flushleft}
\caption{Estimators of $\tau_1$ for several simulated scenarios, when the size of the table is $C=6\cdot 10^5$.}
\label{tab:simulated_example_m}
\end{flushleft}
\end{center}
\end{table}

\begin{table}[h!]
\begin{center}
\begin{tabular}{{|c||ccccccc|}}\hline
{}  & {Zipf $0.2$} & {Zipf $0.5$} & {Zipf $0.8$} & {Zipf $1$}  & {Uniform} & {Dirichlet $0.5$} & {Dirichlet $1$} \\ \hline\hline
True $\tau_1$   & 28976  & 27049  &  21933 &  13794  & 29483  & 15635 & 20281 \\\hline
$\hat{\tau}_1^{L_b}$ & 49619 & 41086  & 28366  & \textbf{ 15980}  & 51607  & \textbf{17952} & 29496\\[0.16cm]
$\hat{\tau}_1^{L_p}$ & 63328 & 53867  & 35646  &  19764  & 64964 & 34447  & 46199\\[0.16cm]
$\hat{\tau}_1^{\Ncr}$ &8076 & 7406   &  5082 & 2729 & 8138 & 6729 & 7371\\[0.16cm]
$\hat{\tau}_1^{\Dcr}$ & 42337 & \textbf{27383}  & 8651 & 3246 & 44104  & 20789 & \textbf{28307}\\[0.16cm]
$\hat{\tau}_1^{B}$ & \textbf{34628}  & 32675  &  \textbf{18977} & 6942  &  \textbf{34772}  & 31235  & 32935\\[0.16cm]
$\hat{\tau}_1^{S}$ & 35957 & 32338  & 15028  & 4196 & 36234  & 29837 & 32804\\\hline
\end{tabular}
\begin{flushleft}
\caption{Estimators of $\tau_1$ for several simulated scenarios, when the size of the table is $C=9\cdot 10^5$.}
\label{tab:simulated_example_l}
\end{flushleft}
\end{center}
\end{table}

From Tables \ref{tab:simulated_example_s}--\ref{tab:simulated_example_l}, we observe that the choice of the smoothing distribution $L$ for  $\hat{\tau}_1^L$, i.e. the Binomial smoothing or the Poisson smoothing, is crucial with respect to the performance of the corresponding estimators. In particular, in all the simulated scenarios the Binomial smoothing displays a better performance than the Poisson smoothing. We also observe that the performance of the estimators strongly depends on the size of the contingency table along with the distributions over the cell's probabilities. However, there is not a clear path indicating in which simulated scenarios nonparametric estimators outperform parametric estimators. In general, we may say that nonparametric estimators have a better performance than parametric estimators when the size of the contingency table is relative small. This confirm a phenomenon already observed in the experimental analysis presented in \citet{Sam(98)} for $\hat{\tau}^{\Dcr}_1$. In general, the nonparametric estimator $\hat{\tau}_1^{L_b}$ has a good performance when the cells probabilities follow the uniform Dirichlet distribution, the uniform distribution and the Zipf with parameter $1$. Also, $\hat{\tau}_1^{L_b}$ appears to be more accurate when the size of the contingency table is relative small.

\clearpage
\appendix
\numberwithin{equation}{section}
\numberwithin{thm}{section}
\numberwithin{lem}{section}
\numberwithin{prp}{section}

\section{Nonparametric estimators of the disclosure risk: proofs}

Here we will prove all the results stated in Section \ref{sec:theory}.
For the sake of simplifying notations, we will simply write $\tau_1$ instead of $\tau_1 (\bm{X},N,M)$, as well as
$\hat{\tau}_1$ (resp. $\hat{\tau}_1^L$) instead of $\hat{\tau}_1 (\bm{X}(N),N)$ (resp. $\hat{\tau}_1^L(\bm{X}(N),N)$).

\subsection{Details for the determination of \eqref{eq:comparingE}} \label{sec:comparingE}

First of all observe that the expected value of $Z_i(\bm{X},N)$ can be easily computed
\begin{equation} \label{eq:EZi}
\begin{split}
\E [Z_i(\bm{X},N)] &= \E \left[ \sum_{j \geq 1} \mathds{1}_{\{Y_j(\textbf{X},N)=i\}}\right] =
\sum_{j \geq 1} \P (Y_j(\textbf{X},N)= i )\\
&=\sum_{j \geq 1}  e^{-np_j} \frac{(np_j)^i}{i!} .
\end{split}
\end{equation}
Using the Taylor series expansion of the exponential function, we get
\begin{align*}
\E [\tau_1] 
&= \E \left[\sum_{j\geq 1} \mathds{1}_{\{Y_j(\textbf{X},N)=1\}}\mathds{1}_{\{Y_j(\textbf{X},N+M)=1\}} \right]\\
&= \sum_{j\geq 1} \P (Y_j (\textbf{X},N) =1) \P (Y_j ( \textbf{X},N+ M) - Y_j (\textbf{X},N)=0  )\\
&= \sum_{j \geq 1} np_j  e^{-n p_j } e^{-\lambda n p_j} = \sum_{j \geq 1}  np_j  e^{-n p_j }
\sum_{i \geq 0} \frac{(-\lambda n p_j)^i}{i!}\\
&= \sum_{i \geq 0} \frac{(-1)^i \lambda^i}{i!} \sum_{j \geq 1} (np_j)^{i+1} e^{-np_j}
 = \sum_{i \geq 0}(-1)^i \lambda^i (i+1) \E [Z_{i+1} (\bm{X},N)]
\end{align*}
where the last equality follows from \eqref{eq:EZi}.

\subsection{Empirical Bayes approach to determine \eqref{eq:estimator<1}} \label{sec:empirical_Bayes}

The estimator $\hat{\tau}_1$ can be derived as the empirical Bayes estimator
of $\E [\hat{\tau}_1]$ in the sense of \cite{Rob(56)}, see also \cite{Mar(89)} for an overview of the empirical Bayes approach. First, it is worth noticing that
the expectation of $\tau_1$ coincides with
\begin{equation}
\label{eq:ET}
\E [\tau_1] = \sum_{j=1}^{+\infty} e^{-(\lambda +1)n p_j} np_j.
\end{equation}
We observe that the statistic $\tau_1$ is a function of the observations only through the frequency counts $Y_j(\bm{X},N)$, which, in our model, are Poisson distributed 
with parameter $np_j$. In order to derive  the nonparametric empirical Bayes estimator of \eqref{eq:ET}, we assume that $p_1, p_2, \ldots$ are independent and distributed according to the empirical cumulative distribution function $G(p)$ of $p_{i_1}, \ldots , p_{i_k}$, corresponding to the $k$ distinct cells arising from the cross classification of the initial sample, namely $G(p):= k^{-1} \sum_{1\leq t \leq k} \mathbbm{1}_{\{p_{i_t}\leq p\}}$.
Consider a cell $j$ containing $x$ individuals out of the initial sample of size $N$, where $x \geq 0$, then by a proper adaptation of \cite[formula (9)]{Rob(56)} to our setting, we find out that
\begin{equation}
\label{eq:Bayes_est}
 \varphi_{n}(x) := \frac{\int e^{-(\lambda+1)np}np e^{-np} \frac{(np)^x}{x!} G(\D p)}{\int e^{-np}\frac{(np)^x}{x!} G(\D p)}
\end{equation}
is the Bayes estimator of the quantity $e^{-(\lambda +1)n p_j} np_j$ appearing in \eqref{eq:ET} for a cell $j$ which contains $x$ individuals out of the initial sample of size $N$.
We now rewrite $\varphi_n(x)$ as follows
\begin{align*}
\varphi_n (x) &=\frac{\int e^{-(\lambda+1)np}np e^{-np} \frac{(np)^x}{x!} G(\D p)}{\int e^{-np}\frac{(np)^x}{x!} G(\D p)}\\
&= \frac{\int \sum_{i \geq 0} \frac{(-(\lambda+1)np)^i}{i!} np e^{-np} \frac{(np)^x}{x!} G(\D p)}{\int e^{-np}\frac{(np)^x}{x!} G(\D p)} \\
& = \frac{\sum_{i \geq 0} \frac{(-(\lambda+1))^i}{i! x!}(x+i+1)! \int \frac{(np)^{x+i+1}}{(x+i+1)!}e^{-np} G(\D p)}{\int e^{-np}\frac{(np)^x}{x!} G(\D p)}\\
& = \frac{\sum_{i \geq 0} \frac{(-(\lambda +1))^i}{i! x!}(x+i+1)! \E [Z_{x+i+1}(\bm{X},N)]}{\E [Z_x(\bm{X},N)]}.
\end{align*}
Then the nonparametric Bayes estimator of $\E [\tau_1]$  may be obtained summing up over all the possible cross classification of the observed cells:
\begin{equation}
\label{eq:Bayes_est}
\begin{split}
&\sum_{x \geq 0 }  Z_{x}(\bm{X},N)   \varphi_{n}(x) \\
&\qquad =
\sum_{x \geq 0} Z_{x}(\bm{X},N)\frac{\sum_{i \geq 0} \frac{(-(\lambda +1))^i}{i! x!}(x+i+1)! \E [Z_{x+i+1}(\bm{X},N)]}{\E [Z_x(\bm{X},N)]}
\end{split}
\end{equation}
The empirical Bayes estimator of $\E [\tau_1]$ coincides with \eqref{eq:Bayes_est} where we replace the expectations $\E [Z_x(\bm{X},N)]$ with their empirical counterparts $Z_{x}(\bm{X},N)$:
\begin{align*}
\hat{\tau}_1 &= \sum_{x \geq 0} Z_x(\bm{X},N) \frac{\sum_{i \geq 0} \frac{(-(\lambda +1))^i}{i! x!}(x+i+1)! Z_{x+i+1}(\bm{X},N)}{Z_{x}(\bm{X},N)}\\
& =  \sum_{x \geq 0} \sum_{i \geq 0} \frac{(-(\lambda +1))^i}{i! x!}(x+i+1)!Z_{x+i+1}(\bm{X},N)\\
& = \sum_{x \geq 0} \sum_{i \geq x} \frac{(-(\lambda +1))^{i-x}}{(i-x)! x!}(i+1)! Z_{i+1}(\bm{X},N)\\
& = \sum_{i \geq 0} (i+1)Z_{i+1}(\bm{X},N)\sum_{x =0}^i \frac{i!}{(i-x)! x!}(-(\lambda +1))^{i-x} \\
& = \sum_{i \geq 0 } (-1)^i \lambda^i (i+1) Z_{i+1}(\bm{X},N) ,
\end{align*}
hence \eqref{eq:estimator<1} now follows.

\subsection{Proof of Theorem \ref{thm:var_lambda<1}} \label{sec:var_lambda<1}

The unbiasedness of $\hat{\tau_1}$ follows from \eqref{eq:comparingE}.
Hence we focus on the proof of the variance bound \eqref{bound_var_T}.
Thanks to the independence of the  random variables $\left\{ Y_j (\bm{X},N) \right\}_{j \geq 1}$, we may
write the variance $\Var (\tau_1-\hat{\tau}_1)$ as
\begin{equation*}
\sum_{j \geq 1} \Var \left(  \sum_{i \geq 0} (-1)^i (i+1) \lambda^i \ind{\{
Y_j (\bm{X},N) = i+1\}}  - \ind{\{ Y_j (\bm{X},N)=1 \}} \ind{\{ Y_j (\bm{X},N+M)=1 \}}  \right).
\end{equation*}
Now the unbiasedness of the estimator implies
\begin{align*}
&\Var (\tau_1-\hat{\tau}_1)\\
& \quad  =
\sum_{j \geq 1} \E\left[ \sum_{i \geq 0} (-1)^i (i+1) \lambda^i \ind{\{
Y_j (\bm{X},N) = i+1\}}  - \ind{\{ Y_j (\bm{X},N)=1 \}} \ind{\{ Y_j (\bm{X},N+M)=1\}}  \right]^2\\
& \quad = \sum_{j \geq 1} \E \left[ \sum_{i \geq 1} a_i \ind{\{
Y_j (\bm{X},N) = i+1\}} +  \ind{\{ Y_j (\bm{X},N)=1 \}} \left(  a_0 -\ind{\{ Y_j (\bm{X},N+M)=1 \}}  \right) \right]^2,
\end{align*}
where we have defined
\[
a_i := (-1)^i (i+1) \lambda^i .
\] 
It is now easy to observe that the events $\{ ( Y_j (\bm{X},N)=i ) \}_{i \geq 1}$ are all disjoint, hence the variance
$\Var (\tau_1-\hat{\tau}_1)$ may be rewritten as
\begin{align*}
 & \sum_{j \geq 1} \E \left[ \sum_{i \geq 1} a_i^2 \ind{\{
Y_j (\bm{X},N) = i+1\}} +  \ind{\{ Y_j (\bm{X},N)=1 \}} \left(  a_0 -\ind{\{ Y_j (\bm{X},N+M)=1 \}}  \right)^2 \right] \\
& \qquad\qquad=  \sum_{j \geq 1} \E \left[ \sum_{i \geq 0} a_i^2 \ind{\{
Y_j (\bm{X},N) = i+1\}} - \ind{\{ Y_j (\bm{X},N)=1 \}} \ind{\{ Y_j (\bm{X},N+M)=1 \}}  \right]
\end{align*}
observing that $a_0=1$. Thus, simple calculations show that we can bound the variance
as follows
\begin{align}
\Var (\tau_1-\hat{\tau}_1) 
& \leq  \max_{j\geq 0} |a_j|^2 \E [ Z_{\bar{1}} (\bm{X},N) ] -  \sum_{j \geq 1} 
e^{-n(\lambda+1)p_j} np_j \nonumber\\
&= \max_{i \geq 0} |a_i|^2 \E [Z_{\bar{1}} (\bm{X},N) ] -\frac{1}{\lambda+1}\E [Z_{1} (\bm{X},N+M) ] . \label{eq:bound_variance_t_max}
\end{align}
To conclude the proof, it remains to show  that  the $a_i$'s have a maximum for $\lambda <1$, which is attained when $i=i^*:= \lfloor (2\lambda-1)/(1-\lambda)\rfloor \vee 0$.  Hence 
the thesis follows by \eqref{eq:bound_variance_t_max}, realizing that
$\max_{i \geq 0} |a_i| = \Psi (\lambda)$.

\subsection{Proof of Theorem \ref{thm:var+bias_L}} \label{sec:teorema_var+bias_L}

First we focus on the determination of the bound \eqref{eq:bias>1}, concerning the bias. Remember the definition 
of both  $\hat{\tau}_1^L$ and $\tau_1$ to write
\begin{align*}
\E [\hat{\tau}_1^L -\tau_1] &= \E \left[  \sum_{i \geq 0} (-1)^{i} (i+1) \lambda^i
\P (L \geq i) Z_{i+1} (\bm{X},N) \right. \\
& \qquad\qquad\qquad \left.- \sum_{j \geq 1} \ind{\{ Y_j (\bm{X},N)=1 \}} \ind{\{  Y_j (\bm{X},N+M) =1\}}\right]\\
&  = - \E \left[ \sum_{i \geq 0} (-1)^{i} (i+1) \lambda^i \P (L \leq i-1)  Z_{i+1} (\bm{X},N) \right]
\end{align*}
where we have observed that non--smoothed estimator $\hat{\tau}_1$ is unbiased.
It is now easy to see that
\begin{align}
\E [\hat{\tau}_1^L -\tau_1]  & =  - \E \left[ \sum_{i \geq 0} (-1)^{i} (i+1) \lambda^i \P (L \leq i-1)  Z_{i+1} (\bm{X},N) \right]\notag\\
& = - \E \left[ \sum_{i \geq 1} (-1)^{i} (i+1) \lambda^i\P (L \leq i-1) 
\sum_{j \geq 1} \ind{\{ Y_{j} (\bm{X},N) =i+1 \}} \right]\notag\\
&= - \sum_{i \geq 1} \sum_{j \geq 1}
 (-1)^{i} (i+1) \lambda^i \P (L \leq i-1) \P ( Y_j (\bm{X}, N) =i+1)\notag\\
&= - \sum_{i \geq 1} \sum_{j \geq 1}
 (-1)^{i} (i+1) \lambda^i \P (L \leq i-1) e^{-np_j} \frac{(np_j)^{i+1}}{(i+1)!}\notag
 \\
& = - \sum_{j \geq 1} e^{-np_j} np_j \sum_{i \geq 1} (-1)^i \frac{(\lambda n p_j)^{i}}{i!} \P (L \leq i-1). \label{eq:bias1_L}
\end{align}
Now we focus on the evaluation of the sum with respect to $i$, for the sake of clarity 
we write $y:= \lambda n p_j$, hence 
\begin{align*}
&\sum_{i \geq 1} (-1)^i \frac{y^{i}}{i!} \P (L \leq i-1) =
\sum_{i =1}^{+\infty} (-1)^i \frac{y^{i}}{i!} \sum_{k=0}^{i-1}\P (L = k)\\
& \qquad\qquad\qquad
= \sum_{k=0}^{+\infty} \P (L=k) \sum_{i=k+1}^{+\infty} \frac{(-y)^i}{i!}
\end{align*}
and remembering the definition of the incomplete gamma function we obtain that
\begin{align*}
\sum_{i \geq 1} (-1)^i \frac{y^{i}}{i!} \P (L \leq i-1) & =
\sum_{k=0}^{+\infty} \P (L=k) \frac{e^{-y}}{k!} \int_0^{-y} \tau^k e^{-\tau}\D \tau
\\
&= -\sum_{k=0}^{+\infty} \P (L=k) \frac{e^{-y}}{k!} \int_0^y (-s)^k e^s \D s\\
& = - e^{-y } \int_0^y e^s \E_L \left[ \frac{(-s)^L }{L!} \right] \D s.
\end{align*}
Putting the previous expression in \eqref{eq:bias1_L} and observing that
$y = \lambda np_j$, \eqref{eq:bias>1} immediately follows.\\
We are ready to bound the variance of the difference between $\tau_1$ and its estimator  $\hat{\tau}_1^L$.
Recalling that the random variables $\left\{ Y_j(\bm{X},N)  \right\}_{j \geq 1}$ are independent, a direct calculation shows that
\begin{align*}
\Var (\hat{\tau}_1^L-\tau_1)
& = 
\Var \left( \sum_{i \geq 0} (-1)^i (i+1)\lambda^i  Z_{i+1} (\bm{X}, N) \P (L \geq i) \right.\\
& \qquad\qquad\qquad \left.-
\sum_{j=1}^{+\infty} \ind{\{  Y_j (\bm{X}, N) =1 \}} \ind{\{ Y_j (\bm{X}, N +M)= 1 \}} \right)\\
&  = \sum_{j=1}^{+\infty} \Var \Big( \sum_{i =0}^{+\infty} (-1)^i (i+1)
\lambda^i \P (L \geq i) \ind{\{ Y_j (\bm{X}, N)=i+1 \}} \\
& \qquad\qquad\qquad - \ind{\{ Y_j (\bm{X}, N)=1 \}} \ind{\{ Y_j (\bm{X}, N+M)= 1 \}} \Big)\\
&  = \sum_{j=1}^{+\infty} \Var \left( \sum_{i =0}^{+\infty} a_i \ind{\{ Y_j (\bm{X},N)=i+1 \}} - \ind{\{ Y_j (\bm{X},N)=1 \}} \ind{\{ Y_j (\bm{X},N+M)= 1 \}} \right),
\end{align*}
having defined 
\[
a_i :=(-1)^i (i+1)\lambda^i \P (L \geq i)
\]
for any $i \geq 0$. Hence the variance of $\hat{\tau}_1^L-\tau_1$ may be upper bounded by the quantity
\begin{align*}
& \sum_{j=1}^{+\infty} \E \left[
\left( \sum_{i =0}^{+\infty} a_i \ind{\{ Y_j (\bm{X},N)=i+1 \}} - \ind{\{Y_j (\bm{X},N)=1 \}} \ind{\{ Y_j (\bm{X},N+M)= 1 \}} \right)^2 \right]\\
& \;=\sum_{j=1}^{+\infty} \E \left[
\left( \sum_{i =1}^{+\infty} a_i \ind{\{ Y_j (\bm{X},N)=i+1 \}} + \ind{\{ Y_j (\bm{X},N)=1 \}} (a_0-\ind{\{ Y_j (\bm{X},N+M)= 1 \}}) \right)^2 \right]\\
&\;=  \sum_{j=1}^{+\infty} \E \left[
\sum_{i =1}^{+\infty} a_i^2 \ind{\{ Y_j (\bm{X},N)=i+1 \}} + \ind{\{Y_j (\bm{X},N)=1 \}} (a_0-\ind{\{ Y_j (\bm{X},N+M)= 1 \}})^2 \right]
\end{align*}
where we have used the incompatibility of the events $\left\{ ( Y_j (\bm{X},N)=i ) \right\}$
for different values of $j$. We can proceed with the upper bound for the variance as follows
\begin{align}
\Var (\hat{\tau}_1^L-\tau_1) &  = \sum_{j=1}^{+\infty} \E \left[
\sum_{i =0}^{+\infty} a_i^2 \ind{\{ Y_j (\bm{X}, N)=i+1 \}} -\ind{\{ Y_j (\bm{X}, N)=1 \}} \ind{\{ Y_j (\bm{X}, N+M)= 1 \}} \right]\notag\\
& \leq \max_{i \geq 0}|a_i|^2 \E [ Z_{\bar{1}} (\bm{X},N)]- \sum_{j=1}^{+\infty}\E \left[ \ind{\{Y_j (\bm{X}, N) =1 \}} \ind{\{ Y_j (\bm{X}, N+M)= 1 \}} \right]\notag\\
& = \max_{i \geq 0}|a_i|^2 \E [Z_{\bar{1}} (\bm{X},N)]- \sum_{j=1}^{+\infty}
e^{-\lambda np_j} e^{-np_j}np_j\notag\\
& = \max_{i \geq 0}|a_i|^2 \E [Z_{\bar{1}} (\bm{X},N)]- \frac{1}{\lambda+1} \E [Z_{1} (\bm{X},N+M)]. \label{eq:var_L_con_max}
\end{align}
We can estimate the maximum value of the $|a_i|$'s as follows
\begin{align*}
\max_{i\geq 0} |a_i| &= \max_{i\geq 0} (i+1) \lambda^i \P (L \geq i) =
\max_{i \geq 0} (i+1)\lambda^i \sum_{k=i}^{+\infty } \P (L =k)\\
&\leq \max_{i \geq 0} \sum_{k=i}^{+\infty }(i+1)\lambda^i  \P (L =k) 
\leq \sum_{k=0}^{+\infty }(k+1)\lambda^k  \P (L =k)  \\
& = \E_L [(L+1)\lambda^L].
\end{align*} 
Hence, replacing  $\max_{i\geq 0} |a_i|$
with $\E_L [(L+1)\lambda^L]$ in \eqref{eq:var_L_con_max}, the upper bound for the variance becomes
\begin{equation}
\label{eq:var_L}
\Var (\hat{\tau}_1^L-\tau_1) \leq (\E_L [(L+1)\lambda^L])^2\E [Z_{\bar{1}} (\bm{X},N)]-
 \frac{\E [Z_{1} (\bm{X},N+M)]}{\lambda+1} .
\end{equation}
Putting together the bound for the variance \eqref{eq:var_L} and for the bias \eqref{eq:bias>1}, the bound on the MSE \eqref{eq:MSE_bound} easily follows.

\subsection{Proof of Proposition \ref{prp:poisson}} \label{sec:MSE_poiss}

Let us now prove the bound \eqref{eq:MSE_poisson} on the MSE, in order to do this, we use Theorem \ref{thm:var+bias_L}, bounding the two terms appearing in \eqref{eq:MSE_bound} separately.\\
To obtain an estimate
 of first term on the r.h.s. of \eqref{eq:MSE_bound},  we note that for any $y >0$ the following holds
\begin{align*}
-e^{-y} \int_0^y e^s \E_L \left[ \frac{(-s)^L }{L!} \right] \D s &= 
-e^{-y} \int_0^y e^s \sum_{k=0}^{+\infty} e^{-\beta} \frac{\beta^k}{k!} \frac{(-s)^k}{k!} \D s \\
& = -e^{-y-\beta} \int_0^y e^s \sum_{k=0}^{+\infty} \frac{(\beta s)^k (-1)^k}{\Gamma (k+1) k!}  \D s
\end{align*}
Recall  that the Bessel polynomial (see \citet{abramo}) is defined as
\[
J_0(z) := \sum_{k=0}^{+\infty} \frac{(-1)^k z^{2k}}{2^{2k} \Gamma (k+1) k!}.
\]
and that $|J_0(z)|\leq 1$, hence we obtain
\begin{align*}
&\left|-e^{-y} \int_0^y e^s \E_L \left[ \frac{(-s)^L }{L!} \right] \D s \right| \\
& \qquad \qquad\leq e^{-(y+\beta)} \int_0^y e^s |J_0 (2\sqrt{s\beta})| \D s
 \leq e^{-\beta} (1-e^{-y}).
\end{align*}
The previous estimate may be applied to bound the first term on the r.h.s. of \eqref{eq:MSE_bound}, with
$y = \lambda n p_j$, indeed 
\begin{equation}\label{eq:bias_poisson}
\begin{split}
&\left|\sum_{j\geq 1} e^{-p_j n(\lambda +1)} p_j n \int_{0}^{\lambda np_j} e^s \E_{L} 
\left[ \frac{(-s)^L}{L!} \right] \D s \right|\\
& \qquad\qquad \leq \sum_{j \geq 1} e^{-np_j} np_j e^{-\beta} (1-e^{-\lambda np_j})
\leq e^{-\beta} \sum_{j=1}^{+\infty} e^{-np_j} np_j \\
& \qquad\qquad = e^{-\beta} \E [Z_1 (\bm{X}, N)] \leq e^{-\beta} \E [N ] = e^{-\beta} n .
\end{split}
\end{equation}
having observed that the maximum number of species with frequency one in a sample of size $N$ is exactly $N$.\\
Second, in order to upper bound the other term on the r.h.s of \eqref{eq:MSE_bound},
we observe that
\begin{align*}
\E_L [(L+1)\lambda^L] & = \sum_{k=0}^{+\infty} e^{-\beta} 
\frac{\beta^k}{k!} \lambda^k (k+1)\\
&= e^{-\beta} \left( \sum_{k=1}^{+\infty} \frac{(\beta \lambda)^k}{(k-1)!}
+\sum_{k=0}^{+\infty} \frac{(\beta \lambda)^k}{k!} \right)\\
& = e^{-\beta} (e^{\beta\lambda}+\beta\lambda e^{\beta \lambda})
= e^{\beta (\lambda-1)} (1+\beta\lambda),
\end{align*} 
hence we get
\begin{equation}
\label{eq:var_poisson}
\begin{split}
&(\E_L [(L+1)\lambda^L])^2 \E [Z_{\bar{1}} (\bm{X},N)] -\frac{1}{\lambda+1} \E [Z_1 (\bm{X}, N+M)]\\
&\qquad\qquad\qquad\qquad\qquad\qquad\qquad\qquad \qquad\leq ne^{2\beta (\lambda-1)} (1+\beta\lambda)^2 . 
\end{split}
\end{equation}
Using \eqref{eq:bias_poisson} and \eqref{eq:var_poisson}, one can now estimate the MSE \eqref{eq:MSE_bound} in the Poisson case and \eqref{eq:MSE_poisson} follows.\\

Thanks to \eqref{eq:MSE_poisson} just derived, the normalized mean 
square error can be bounded from above by
\begin{equation*}
\Ecr_{n, \lambda} (\hat{\tau}^L_1) \leq e^{-2\beta} +\frac{e^{2\beta (\lambda-1)} (1+\beta \lambda)^2}{n}
\end{equation*}
using the exponential inequality $1+x\leq e^x$ we get
\begin{equation} \label{eq:NMSE_poisson_for_minimizing}
\Ecr_{n, \lambda} (\hat{\tau}^L_1) \leq e^{-2\beta} +\frac{e^{2\beta (2\lambda-1)}}{n}.
\end{equation}
The r.h.s. of \eqref{eq:NMSE_poisson_for_minimizing} is minimized when $\beta$
equals $\frac{1}{4\lambda}\log \left(  \frac{n}{2\lambda-1}\right)$,
it is easy to observe that \eqref{eq:NMSE_poisson_for_minimizing} becomes
\begin{equation} 
\Ecr_{n, \lambda} (\hat{\tau}^L_1) \leq 
\frac{1}{n^{1/(2\lambda)}}\cdot \frac{2\lambda}{(2\lambda-1)^{1-1/(2\lambda)}}
\end{equation}
hence the second bound \eqref{eq:NMSE_poisson} follows provided that
\[
A(\lambda) :=.\frac{2\lambda}{(2\lambda-1)^{1-1/(2\lambda)}} .
\]
We are now ready to prove the limit of predictability in the Poisson case, indeed thanks to \eqref{eq:NMSE_poisson} we have
\[
\Ecr_{n, \lambda} (\hat{\tau}^L_1) \leq  \frac{A}{n^{1/(2\lambda)}},
\]
besides observe that the inequality
\[
\frac{A}{n^{1/(2\lambda)}} \leq \delta
\]
is satisfied iff 
\[
\lambda \leq \frac{\log(n)}{2 \log (A/\delta)}=:\lambda^*.
\]
As a consequence the maximum value of $\lambda$ for which the inequality
$\Ecr_{n, \lambda} (\hat{\tau}^L_1) \leq \delta$ is satisfied, is bigger or equal than $\lambda^*$, in other words
\[
\max \left\{ \lambda : \; \Ecr_{n, \lambda} (\hat{\tau}^L_1)\leq \delta \right\} \geq \frac{\log(n)}{2 \log (A/\delta)}.
\]
The thesis follows by taking the limit of the previous inequality as $n \to +\infty$.

\section{Proofs of auxiliary results for the lower bound}
\label{sec:proofs-auxil-results}

\subsection{Proof of Lemma \ref{lem:3}}
\label{sec:proof-lem:3}

First, it is obvious that
\begin{equation*}
  \Ecr(\lambda , n)  \leq \inf_{\hat{\rho}} \sup_{P\in
    \Pcr} n^{-2}\E_P^n[(\tau_1(\bm{X},N,M) -
  \hat{\rho}(\bm{Y}(\bm{X},N)))^2].
\end{equation*}
We now prove that the previous is indeed an inequality by deriving a lower bound that essentially  matches. Let $n > 0$ be fixed. By definition, for every $\varepsilon > 0$
there exists an estimator $\hat{\rho}_1$ such that
\begin{align}
  \notag
  \Ecr(\lambda , n) %
  &\geq \sup_{P\in \Pcr}n^{-2}\E_P^n[(\tau_1(\bm{X},N,M) -
    \hat{\rho}_1(\bm{X}(N),N))^2] - \varepsilon\\
  \notag
  &=\sup_{P\in \Pcr}n^{-2}\E_P^n[ \E_P^n[(\tau_1(\bm{X},N,M)\\
  & \qquad\qquad\qquad -
    \hat{\rho}_1(\bm{X}(N),N))^2 \mid \bm{Y}(\bm{X},N),
    \bm{Y}(\bm{X},N+M)]] - \varepsilon \notag\\
  \label{eq:5}
  &\geq \sup_{P\in \Pcr}n^{-2}\E_P^n[(\tau_1(\bm{X},N,M) -
    \E_P^n[\hat{\rho}_1(\bm{X}(N),N)\mid \bm{Y}(\bm{X},N)])^2] -
    \varepsilon
\end{align}
where the last line follows by Jensen's inequality and by observing that
\begin{gather*}
  \E_P^n[\tau_1(\bm{X},N,M) \mid \bm{Y}(\bm{X},N),
  \bm{Y}(\bm{X},N + M)]%
  = \tau_1(\bm{X},N,M),\quad\mathrm{and},\\
  \E_P^n[\hat{\rho}_1(\bm{X}(N),N) \mid \bm{Y}(\bm{X},N),
  \bm{Y}(\bm{X},N + M)] = \E_P^n[\hat{\rho}_1(\bm{X}(N),N)
  \mid \bm{Y}(\bm{X},N)].
\end{gather*}
To see that the last equation is true, remark that
$\bm{Y}(\bm{X},N+M) - \bm{Y}(\bm{X},N)$ is independent of
$\bm{Y}(\bm{X},N)$ and depends only on $(X_{N+1},\dots,X_{N+M})$. Now we
claim that $\hat{\rho}_1$ can be chosen such that for any $k\in
\Z_+$ and any
permutation $\sigma_k(\bm{X}(k))$ of the data, it holds
$\hat{\rho}_1(\bm{X}(k),k) = \hat{\rho}_1(\sigma_k(\bm{X}(k)),k)$. We
delay the proof of the claim to later. Now assume the claim is true. Given $k$
and $\bm{Y}(\bm{X},k)$, we can construct the functional
\begin{equation*}
  G(\bm{Y}(\bm{X},k),k)%
  := ( \underbrace{1, \ldots, 1}_{\times
    Y_1(\bm{X},k)},\underbrace{2,\ldots , 2}_{\times  Y_2(\bm{X},k)},\dots).
\end{equation*}
Since $\hat{\rho}_1$ is invariant under permutations of the data, we have for
any $P \in \Pcr$,
\begin{multline*}
  \E_P^n[\hat{\rho}_1(\bm{X}(N),N) \mid \bm{Y}(\bm{X},N)]\\%
  \begin{aligned}
    &= \E_P^n\big[\E_P^n[\hat{\rho}_1(\bm{X}(N),N) \mid
    \bm{Y}(\bm{X},N),N] \mid \bm{Y}(\bm{X},N) \big]\\
    &= \E_P^n\big[ \E_P^n[\hat{\rho}_1(G(\bm{Y}(\bm{X},N),N),N) \mid
    \bm{Y}(\bm{X},N),N] \mid \bm{Y}(\bm{X},N) \big]\\%
    &= \E_P^n[\hat{\rho}_1(G(\bm{Y}(\bm{X},N),N),N) \mid
    \bm{Y}(\bm{X},N) ]\\
    &= \hat{\rho}_1(G(\bm{Y}(\bm{X},N),N),N).
  \end{aligned}
\end{multline*}
The last line follows because $N = \sum_{j\geq 1}Y_j(\bm{X},N)$, and
hence $N$ is completely determined by $\bm{Y}(\bm{X},N)$. Therefore, we
have proved that the conditional expected value of $\hat{\rho}_1(\bm{X}(N),N)$, given 
$\bm{Y}(\bm{X},N)$ does not depend on $P$. Thus, \eqref{eq:5} implies,
\begin{align*}
  \Ecr(\lambda , n) %
  &\geq \sup_{P\in \Pcr} n^{-2}\E_P^n[ (\tau_1(\bm{X},N,M) -
    \hat{\rho}_1(G(\bm{Y}(\bm{X},N),N),N))^2 ] - \varepsilon\\
  &\geq \inf_{\hat{\rho}} \sup_{P\in \Pcr} n^{-2}\E_P^n[
    (\tau_1(\bm{X},N,M) - \hat{\rho}(\bm{Y}(\bm{X},N)))^2 ] -
    \varepsilon.
\end{align*}
Since the previous is true for all $\varepsilon > 0$, the conclusion follows.\\

We now prove the claim we have used in the previous argument, i.e. that
\textit{$\hat{\rho}_1$ can be chosen such for any $k\in \Z_+$ and any
permutation $\sigma_k(\bm{X}(k))$ of the data, it holds
$\hat{\rho}_1(\bm{X}(k),k) = \hat{\rho}_1(\sigma_k(\bm{X}(k)),k)$.}
 When $k = 0$, then the claim is trivial, hence we assume
without loss of generality that $k \in \N$. We will prove that for any
estimator $\hat{\rho}_1$, there is a symmetric estimator $\hat{t}_1$ with a risk
no more than the risk of $\hat{\rho}_1$. Let $\hat{\rho}_1$ be
arbitrary. Construct $\hat{t}_1$ such that for any $k\in \N$
\begin{equation*}
  \hat{t}_1(\bm{X}(k),k)%
  :=
  \frac{1}{|\{\sigma_k\}|}\sum_{\{\sigma_k\}}\hat{\rho}_1(\sigma_k(\bm{X}(k)),k).
\end{equation*}
Clearly $\hat{t}_1$ has the desired invariance property under
permutations. Moreover, by Jensen's inequality,%
\begin{multline*}
  \E_P^n[(\tau_1(\bm{X},N,M) - \hat{t}_1(\bm{X}(N),N) )^2]\\
  \begin{aligned}
    &= \E_P^n\Big[\E_P^n\Big[\Big(\frac{1}{|\{\sigma_N\}|}\sum_{\{\sigma_N\}}
    (\tau_1(\bm{X},N,M) -
    \hat{\rho}_1(\sigma_N(\bm{X}(N) ),N) \Big)^2 \mid N \Big]\Big]\\
    &\leq \E_P^n \Big[ \E_P^n\Big[ \frac{1}{|\{\sigma_N\}|}\sum_{\{\sigma_N\}}
    (\tau_1(\bm{X},N,M) - \hat{\rho}_1(\sigma_N(\bm{X}(N)),N) )^2 \mid N
    \Big] \Big]
  \end{aligned}
\end{multline*}

Now remark that for all $(k,k') \in \Z_+^2$ the map
$\bm{X} \mapsto \tau_1(\bm{X},k,k')$ is invariant under any permutations
of the $k$ first entries of $\bm{X}$. Moreover, $\bm{X}$ is an i.i.d.
vector, then the last display implies that
\begin{multline*}
  \E_P^n[(\tau_1(\bm{X},N,M) - \hat{t}_1(\bm{X}(N),N) )^2]\\
  \begin{aligned}
    &= \E_P^n \Big[\E_P^n\Big[\frac{1}{|\{\sigma_N\}|}\sum_{\{\sigma_N\}}
    (\tau_1(\bm{X},N,M)
    - \hat{\rho}_1(\bm{X}(N),N) )^2 \mid N\Big]\Big]\\
    &= \E_P^n[(\tau_1(\bm{X},N,M) - \hat{\rho}_1(\bm{X}(N),N) )^2].
  \end{aligned}
\end{multline*}
The conclusion follows by taking the supremum over $P \in \Pcr$ both sides of
the last display.

\subsection{Proof of Proposition \ref{pro:6}}
\label{sec:proof-pro:6}

For any $P \in \Pcr'$ we write $\tilde{P}(\cdot) := P(\cdot) / P(\N)$, so that
$\tilde{P} \in \Pcr$ is a probability measure. We write
$\tilde{p}_j := p_j / P(\N)$, $j \in \{1,\dots,S\}$. Furthermore we let
$m(P) := n \sum_{j=1}^S p_j$. Then since $\bm{Y}$ is a vector of independent
Poisson random variables, is clear that for any $P \in \Pcr'$
\begin{equation}
   \label{eq:1}
  \E_{\tilde{P}}^n[(\bar{\tau}_1(\tilde{P},n) - \hat{\rho}(\bm{Y}))^2]%
  = \E_P^{m(P)}[(\bar{\tau}_1(\tilde{P},n) - \hat{\rho}(\bm{Y}) )^2].
\end{equation}
We now choose $\hat{\tau}$ to be an estimator satisfying for some $\zeta > 0$
\begin{equation*}
  \sup_{P\in \Pcr'}\E_P^{m(P)}[(\bar{\tau}_1(\tilde{P},n) -
  \hat{\tau}(\bm{Y}))^2]%
  \leq \inf_{\hat{\rho}}\sup_{P\in \Pcr'}\E_P^{m(P)}[(\bar{\tau}_1(\tilde{P},n) -
  \hat{\rho}(\bm{Y}))^2] + \zeta.
\end{equation*}
This is always possible for any $\zeta > 0$. Furthermore remark that
$m(P) \leq (1+c_1\varepsilon/\xi)n = n'$, so that $m(P)/n' \leq 1$ always when
$P \in \Pcr'$. Let $P \in \Pcr'$ be fixed, and let $\bm{W}= (W_1,W_2,\dots)$
such that conditional on $\bm{Y}$, the random variables $W_j$ are independent
binomial random variables with parameters $(Y_j, m(P)/n')$. Then define
$\tilde{\tau}(\bm{Y}) \coloneqq \E[\hat{\tau}(\bm{W})] \mid \bm{Y}]$. By
Jensen's inequality,
\begin{align*}
  \E_P^{n'}[(\bar{\tau}_1(\tilde{P},n) - \tilde{\tau}(\bm{Y}))^2]%
  &= \E_P^{n'}[( \E[\bar{\tau}_1(\tilde{P},n) - \hat{\tau}(\bm{W}) \mid \bm{Y}  ] )^2]\\
  &\leq \E_P^{n'}[\E[(\bar{\tau}_1(\tilde{P},n) - \hat{\tau}(\bm{W}))^2 \mid
    \bm{Y}]]\\
  &= \E_P^{m(P)}[(\bar{\tau}_1(\tilde{P},n) - \hat{\tau}(\bm{Y}))^2]\\
  &\leq \inf_{\hat{\rho}}\sup_{P\in \Pcr'}\E_P^{m(P)}[(\bar{\tau}_1(\tilde{P},n) -
  \hat{\rho}(\bm{Y}))^2] + \zeta.
\end{align*}
Taking the supremum over $P \in \Pcr'$ on the lhs of the last display, and using
that the infinimum over $\hat{\rho}$ will be always smaller than the value at
$\tilde{\tau}$, we find using \eqref{eq:1} that
\begin{align*}
  \inf_{\hat{\rho}}\sup_{P\in \Pcr}\E_P^n[(\bar{\tau}_1(P,n) -
  \hat{\rho}(\bm{Y}))^2]%
  &= \inf_{\hat{\rho}}\sup_{P\in \Pcr'}\E_{\tilde{P}}^n[(\bar{\tau}_1(\tilde{P},n) -
    \hat{\rho}(\bm{Y}))^2]\\%
  &= \inf_{\hat{\rho}}\sup_{P\in \Pcr'}\E_P^{m(P)}[(\bar{\tau}_1(\tilde{P},n) -
    \hat{\rho}(\bm{Y}))^2]\\
  &\geq \inf_{\hat{\rho}}\sup_{P\in \Pcr'}\E_P^{n'}[(\bar{\tau}_1(\tilde{P},n) -
  \hat{\rho}(\bm{Y}))^2] - \zeta
\end{align*}
Since the previous is true for all $\zeta > 0$, we indeed have proven
\begin{equation}
  \label{eq:12}
  \inf_{\hat{\rho}}\sup_{P\in \Pcr}\E_P^n[(\bar{\tau}_1(P,n) -
  \hat{\rho}(\bm{Y}))^2]%
  \geq \inf_{\hat{\rho}}\sup_{P\in \Pcr'}\E_P^{n'}[(\bar{\tau}_1(\tilde{P},n) -
  \hat{\rho}(\bm{Y}))^2].
\end{equation}

\par To finish the proof of the proposition, we will now show that
$\bar{\tau}_1(\tilde{P},n)$ in \eqref{eq:12} can be traded for $\bar{\tau}_1(P,n)$
at small cost. Remark that by Young's inequality, for any $P \in \Pcr'$ and any
$\hat{\rho}$,
\begin{equation}
  \label{eq:13}
  \begin{split}
&  \E_P^{n'}[(\bar{\tau}_1(\tilde{P},n) - \hat{\rho}(\bm{Y}))^2 ]\\
 &\qquad\qquad \geq \frac{1}{2} \E_P^{n'}[(\bar{\tau}_1(P,n) - \hat{\rho}(\bm{Y}))^2 ]
  - ( \bar{\tau}_1(P,n) - \bar{\tau}_1(\tilde{P},n))^2,
  \end{split}
\end{equation}
with
\begin{equation*}
\begin{split}
&\bar{\tau}_1(P,n) - \bar{\tau}_1(\tilde{P},n) \\
  & \qquad=
  - n\sum_{j=1}^S (\tilde{p}_j- p_j)e^{-(1 + \lambda)np_j}
  + n \sum_{j=1}^S \tilde{p}_je^{-(1 + \lambda) n p_j}\Big\{1 -
  e^{n(1+\lambda)(p_j - \tilde{p}_j)} \Big\}.
  \end{split}
\end{equation*}
Thanks to a Taylor expansion of the term within the  brackets in the last display for
$p_j - \tilde{p}_j$ near to $0$, we find that for $n$ large enough,
\begin{align*}
  |\bar{\tau}_1(P,n) - \bar{\tau}_1(\tilde{P},n)|
  &\leq c_1n \varepsilon/\xi
    + (c_1 n \varepsilon/\xi) n(1 + \lambda)\sum_{j=1}^S \tilde{p}_j^2\\%
  &\leq (c_1n \varepsilon/\xi) + (c_1n \varepsilon/\xi) n(1+\lambda)
    \max_{j=1,\dots,S} \tilde{p}_j\\
  &= c_1n \varepsilon/\xi + (c_1n \varepsilon/\xi) n(1+\lambda) (1+
    O(c_1\varepsilon/\xi)) \max_{j=1,\dots,S} p_j\\
  &\leq 3 c_1n \varepsilon.
\end{align*}
This estimate combined with \eqref{eq:12}, \eqref{eq:13} and \eqref{eq:11} completes the proof for
the first inequality of the proposition. The second inequality simply follows
from the first by an application of Markov's inequality.

\subsection{Proof of Proposition \ref{pro:2}}
\label{sec:proof-pro:2}

By a simple application of Bernstein's inequality, we get that
\begin{align*}
  \P\big( Q(\bm{U}) \notin \Pcr' \big)
  =\P\big(| \textstyle\sum_{j=1}^SU_j - 1| > c_1\varepsilon/\xi \big)\leq 2\exp\Big\{- \frac{1}{2} \frac{c_1^2 \varepsilon^2/ \xi^2}{\xi S^{-1}
    + \frac{1}{3} S^{-1}c_1\varepsilon } \Big\}.
\end{align*}
Then the conclusion follows from simple algebraic manipulations.

\subsection{Proof of Proposition \ref{pro:7}}
\label{sec:proof-pro:7}

By definition, we have that
\begin{equation*}
  \bar{\tau}_1(Q(\bm{U}),n)%
  = n%
  \sum_{j=1}^S U_j e^{-n(1+\lambda)U_j}.
\end{equation*}
Whence, $\bar{\tau}_1(Q(\bm{U}))$ is a sum of i.i.d. random variables taking
values in $[0, n \xi S^{-1}]$. By Hoedffding's inequality,
\begin{align*}
  \P\big( |\bar{\tau}_1(Q(\bm{U}),n) - \E[\bar{\tau}_1(Q(\bm{U}),n)]|
  > n \varepsilon/2 \big)%
  &\leq 2 \exp\Big\{ -\frac{S\varepsilon^2}{2\xi} \Big\}
\end{align*}
The conclusion follows from simple algebraic manipulations.

\subsection{Proof of Proposition \ref{pro:8}}
\label{sec:proof-pro:8}

\par We first consider the case where $1+\lambda \leq \log n$. Remark that in
that case we have,
\begin{equation}
  \label{eq:26}
  \frac{n' \xi}{S}%
  = \frac{n(1 + c_1\varepsilon/\xi)(2c_0/e)(1+\lambda)\log(n)}{n(1+\lambda)}%
  \leq \frac{(3 c_0/e)\log(n)}{1+\lambda},
\end{equation}
where the last inequality is true for $n$ large enough. Using  
\cite[Lemma~6]{wu2015chebyshev}, and because $0 \leq U,V \leq \xi S^{-1}$ almost surely,
we find that,
\begin{align*}
 &S\TV\Big(\E[{\rm Poiss}(n' U)], \E[{\rm Poiss}(n' V)] \Big)
  \\
  & \qquad\qquad\leq \frac{S}{(L+2)!}\Big(\frac{n' \xi}{2S}\Big)^{L+2}\Big(2 +
    2^{n'\xi/(2S) - L} + 2^{n'\xi/(2 \log(2) S) - L} \Big)\\
  &\qquad\qquad= \frac{2S(1+o(1))}{(L+2)!}\Big(\frac{n' \xi}{2S}\Big)^{L+2},
\end{align*}
where the last line is a consequence of the definition of $L$ and \eqref{eq:26}. Indeed  we
always have%
\begin{equation*}
  \frac{n'\xi}{2S}%
  < \frac{n'\xi}{2S \log(2)}%
  \leq \frac{3c_0}{2e \log 2} \log(n)%
  < 0.8 c_0\log(n) < L,
\end{equation*}
where the last inequality is again also true at least for $n$ large enough,
because $L = c_0 A(\lambda,n) > c_0(1 + o(1)) \log(n)$, for any choice of
$c_0 > 0$ because $a \log a > 0 \Rightarrow a > 1$. \\
Now, observing that 
$(L+2)! > L^2 L!$, we obtain the upper bound
\begin{align*}
 & S\TV\Big(\E[{\rm Poiss}(n' U)], \E[{\rm Poiss}(n' V)] \Big)\\
  &\qquad\qquad\leq \frac{2S(1+o(1))}{L^2}\Big( \frac{n\xi}{2S} \Big)^2%
    \frac{1}{L!}\Big( \frac{n\xi}{2S} \Big)^L\Big( \frac{n'}{n}\Big)^{L+2}\\
  &\qquad\qquad\leq \frac{2eS(1+o(1))}{L^2}\Big( \frac{c_0\log(n)}{e} \Big)^2%
    \frac{1}{L!}\Big( \frac{c_0\log(n)}{e} \Big)^L,
\end{align*}
where the last line follows because $L+2 \leq c_2 \xi$, and hence
$(1+c_1\varepsilon/\xi)^{L+2}\leq \exp\{c_1c_2 \varepsilon\} \leq e$, by the
assumption on $\varepsilon$. Then, by Stirling's formula, whenever $n
\rightarrow \infty$ (and hence $L$),
\begin{align*}
  S\TV\Big(\E[{\rm Poiss}(n' U)], \E[{\rm Poiss}(n' V)] \Big)%
  &\leq \frac{2(1+o(1))}{\sqrt{2\pi c_0}e A(\lambda,n)^{5/2}}%
  \frac{n(1+\lambda)}{\sqrt{\log(n)}}%
  A(\lambda,n)^{-L}.
\end{align*}
The conclusion then follows for $c_3 > 0$ large enough by the definition of
$A(\lambda,n)$, and because when $1 + \lambda \leq \log(n)$ it holds
$A(\lambda,n) = \omega(1+o(1))$ with $\omega$ solution to $\omega \log \omega =
c_0^{-1}$, hence $c_3$ can be chosen to depends only on $c_0$.

\par We now consider the case where $1+\lambda > \log(n)$. Under this constraint
$\xi = (2c_0/e)\log^2(n)$, and proceeding as for \eqref{eq:26} we find that
$n'\xi/S \leq (3c_0/e)$ as long as $n$ gets large enough. Whence, whenever $n$
is large enough we certainly have $n'\xi/S = o(L)$, and still by 
\cite[Lemma~6]{wu2015chebyshev}, and along similar lines as in the previous
paragraph, we get
\begin{align*}
 &S\TV\Big(\E[{\rm Poiss}(n' U)], \E[{\rm Poiss}(n' V)] \Big)\\%
  &\qquad\qquad\leq \frac{2S(1+o(1))}{(L+2)!}\Big(\frac{n \xi}{2S}\Big)^{L+2}%
    \Big(\frac{n'}{n}\Big)^{L+2}\\
  &\qquad\qquad\leq \frac{2e(1+o(1))}{(L+2)!} n(1+\lambda) \Big( \frac{c_0
    \log^2(n)}{e(1+\lambda)} \Big)^{L+2}\\
  &\qquad\qquad\leq \frac{2c_0^2(1+o(1))}{e \sqrt{2\pi}} \frac{n \log^4(n)}{1+\lambda}%
    \frac{1}{L^{5/2}}\Big(\frac{c_0 \log^2(n)}{(1+\lambda) L} \Big)^L\\
  &\qquad\qquad\leq \frac{2c_0^2(1+o(1))}{e \sqrt{2\pi}}%
    \frac{n \log^3(n)}{L^{5/2}} \Big(\frac{c_0 \log(n)}{ L} \Big)^L.
\end{align*}
but $L \geq 2c_0\log(n)$, so that the previous bound always goes to zero
when $n \rightarrow \infty$, and hence gets smaller than $1/2$ for $n$ large
enough.

\subsection{Proof of Proposition \ref{pro:1}}
\label{sec:proof-pro:1}

We define the function $\varphi : \R_{+} \rightarrow \R$ such that
$\varphi(x) = x \log(x)$. When $1 + \lambda \leq \log(n)$, it is clear that
$A(\lambda,n)$ converges to the solution of $\varphi(x) = c_0^{-1} = e$, hence
$A(\lambda,n) \rightarrow e$, which proves the first claim.

\par For the second claim, let define,
\begin{equation*}
  \Delta_n%
  := e\frac{\log(1+\lambda) - (1/2)\log\log(n) + \log(c_3)}{\log(n)}.
\end{equation*}
For $n$ large enough such that $\Delta_n > -1$, it is clear than
$A(\lambda,n) \geq 0$. Furthermore, by a Taylor expansion of $\varphi$ near
$x = e$, we find that there is a $\bar{x}$ in the line segment
between $A(\lambda,n)$ and $e$,
\begin{align*}
  \varphi(A(\lambda,n))%
  &= \varphi(e) + \varphi'(e)(A(\lambda,n) - e) +
    \frac{\varphi''(\bar{x})}{2}(A(\lambda,n) - e)^2\\
  &\geq \varphi(e) + \varphi'(e)(A(\lambda,n) - e),
\end{align*}
because $\varphi''(x) = 1/x > 0$ whenever $x > 0$. Since
$\varphi(A(\lambda,n)) - \varphi(e) = \Delta_n$, $\varphi(e) = e$, and
$\varphi'(e) = 2$, we deduce that for those $n$ large,
\begin{equation*}
  0 \leq A(\lambda,n) \leq e + \Delta_n/2.
\end{equation*}
Therefore,
\begin{align*}
  e^{-1}A(\lambda,n)^2 \log(n)%
  &\leq e \log(n)%
    + \Delta_n \log(n)%
    + \frac{\Delta_n^2 \log(n)}{4e}\\
  &= e \log(n)%
    + e\log \frac{c_3(1+\lambda)}{\sqrt{\log(n)}}%
    +o(1).
\end{align*}
This concludes the proof.

\section{Existence of random variables}
\label{sec:exist-rand-vari-1}

Here we prove the existence of the random variables $U$ and $V$ which have been used to construct the prior 
for determining the  minimax lower bound in Section \ref{sec:proof-lower-bound}. More precisely, we prove the
following theorem.
\begin{thm}
  \label{thm:4}
  Let $S,L\in \N$ and $\xi>0$ chosen as in Section \ref{sec:proof-lower-bound}. Then there exist two random variables
$U$ and $V$ taking
values in $[0,\xi S^{-1}]$ such that when $n$ is large enough,
\begin{gather*}
  \E[U^k] = \E[V^k]\quad \forall k \in \{0,\dots,L+1\},\\
  \E[U] = \E[V] = S^{-1},\quad%
  \Var(U) \leq \xi S^{-2},\quad \Var(V) \leq \xi S^{-2},\\
  \E[U e^{-n(1+\lambda)U}] \geq \E[V e^{-n(1+\lambda)V}]%
  + S^{-1} K\min\big\{1,\, \sqrt{\xi/L^2}\exp(- L^2/\xi) \big\}.
\end{gather*}
\end{thm}
The proof of Theorem \ref{thm:4} follows the guidelines used in the papers
\citet{wu2015chebyshev,wu2016minimax}, relating the problem of the existence of
the random variables to the problem of finding the best polynomial approximation
to some function.\\

For $a,b \in \R$, we let $\cont[a,b]$ denote the space of continuous
functions on $[a,b]$, and for any $L \in \Z_+$ we let
$\pol_L[a,b] \subset \cont[a,b]$ denote the space of polynomials of degree no
more than $L$ on $[a,b]$. For any $f \in \cont[a,b]$, the best polynomial (of
degree at most $L$) approximation to $f$ is defined as
\begin{equation*}
  E_L(f, [a,b]) :=%
  \inf\{ \sup\{|f(x) - q(x)| : \;  x \in [a,b]\} : \;  q \in
    \pol_L[a,b]\}.
\end{equation*}

\par For the sake of simplicity, we define $B :=
n(1+\lambda)\xi/(2S)$. Remark that $B \asymp \xi/2$, but is not necessarily
equal to it because $S$ is integer. We define
$g : [\xi^{-1},1] \rightarrow \R_{+}$ such that
$g(x) := \exp\{- 2B x\}$. It is a classical result that for any
$L \in \N$ we can find random variables $X$ and $Y$ taking values in
$[\xi^{-1},1]$ and such that
\begin{gather*}
  \E[X^k] = \E[Y^k],\qquad k=0,\dots, L,\\
  \E[g(X)] = \E[g(Y)] + E_L(g, [\xi^{-1},1]).
\end{gather*}
The proof of the existence of such random variables can be found for instance in
\citet{wu2016minimax,wu2015chebyshev} for a constructive argument, or for
instance in \citet{Lep(99)} using the Hahn-Banach theorem and a
duality argument.

\par We now assume that we have at our disposal the random variables $X$ and $Y$
of the previous paragraph, and we write $P_X$ and $P_Y$ their distributions. The
construction of the random variables $U$ and $V$ is done using the trick
introduced in \citet[Lemma~4]{wu2016minimax}. Namely, we let $U$ and $V$ having
respective distributions on $[0, \xi S^{-1}]$
\begin{align*}
  P_U(\D x)%
  &:= \big(1 - \E[(\xi X)^{-1}]\big) \delta_0 + (S x)^{-1}P_{\xi
    X/S}(\D x),\\
  P_V(\D x)%
  &:= \big(1 - \E[(\xi Y)^{-1}]\big) \delta_0 + (S x)^{-1}P_{\xi
    Y/S}(\D x).
\end{align*}
Because $X,Y \geq \xi^{-1}$ almost-surely, then $\E[(\xi X)^{-1}] \leq 1$ and
$\E[(\xi Y)^{-1})] \leq 1$. Indeed from \citet[Lemma~4]{wu2016minimax}, $P_U$
and $P_V$ are proper probability distributions on $[0, \xi S^{-1}]$ satisfying
\begin{align*}
  &\E[U] = \E[V] = 1/S,\qquad%
  \E[U^k] = \E[V^k],\qquad k=0,\dots, L+1,\\
  &\E[U \exp\{- n(1+\lambda) U\} ] = \E[V \exp\{-n(1+\lambda)V\}] +
  S^{-1}E_{L}(g, [\xi^{-1},1]).
\end{align*}
Furthermore, it is clear that,
\begin{equation*}
  \E[U^2]%
  = \frac{1}{S}\int x\, P_{\xi x/S}(\D x)%
  = \frac{\xi \E[X]}{S^2}%
  \leq \frac{\xi}{S^2}.
\end{equation*}
Hence $\Var(U) \leq \xi / S^2$. It is obvious that we also have
$\Var(V) \leq \xi/S^2$. Thus, the proof of the theorem is finished by
obtaining a lower bound on the best polynomial approximation
$E_L(g,[\xi^{-1},1])$, which is done in the next section (in particular see
Theorem \ref{thm:5} and the \eqref{eq:25} just after).

\section{Approximation theory}
\label{sec:approximation-theory}

\subsection{Statement of the main result}
\label{sec:stat-main-result-2}

We wish to find the best polynomial approximation (see
Section \ref{sec:exist-rand-vari-1} for definition) to the function
$g : [\xi^{-1},1] \rightarrow \R_{+}$ such that
$g(x) := \exp\{- 2B x\}$ on $[\xi^{-1},1]$, with $\xi$ defined in
\eqref{eq:21} and $B$ satisfying
\begin{equation}
  \label{eq:24}
  \frac{\xi/2}{1 + \frac{1}{n(1+\lambda)}} \leq B \leq \xi/2.
\end{equation}
The whole section will be dedicated to the proof of the following theorem.
\begin{thm}
  \label{thm:5}
  For every $\zeta > 0$, there exists a constant $K > 0$ such that
  as $n \rightarrow \infty$,
  \begin{equation*}
    E_L(g, [\xi^{-1},1]) \geq%
    K(1+o(1))\cdot
    \begin{cases}
      1 &\mathrm{if}\ L \leq \sqrt{\xi/2},\\
      \frac{\sqrt{\xi}\exp\{-L^2/\xi\} }{L(1+(2L/\xi)^2)^{1/4}}
      &\mathrm{if}\ \sqrt{\xi/2} < L < \zeta \xi.
    \end{cases}
  \end{equation*}
\end{thm}

We deduce from the previous theorem that there exists a universal constant $K >
0$ such that for $n$ large enough,
\begin{equation}
  \label{eq:25}
  E_L(g, [\xi^{-1},1])%
  \geq K\min\big\{ 1,\, \sqrt{\xi/L^2}\exp(-L^2/\xi) \big\}.
\end{equation}

\subsection{Proof of Theorem \ref{thm:5}}
\label{sec:proof-thm:5}

\par Let $\sigma : [-1,1] \rightarrow [\xi^{-1},1]$ be such that
$\sigma(x) := (1-\xi^{-1})(x+1)/2 + \xi^{-1}$. Notice that $\sigma$ is
bijective. By translating and rescaling, we claim that
$E_L(g, [\xi^{-1},1]) = E_L(g \circ \sigma, [-1,1])$. To see that this is true,
remark that for all $p \in \pol_L[-1,1]$ we have
$\|g\circ \sigma - p\|_{\infty} = \|g - p \circ \sigma^{-1}\|_{\infty} \geq
E_L(g,[\xi^{-1},1])$. This shows that
$E_L(g \circ \sigma,[-1,1]) \geq E_L(g,[\xi^{-1},1])$. The same steps using
$\sigma^{-1}$ show that $E_L(g \circ \sigma,[-1,1]) \leq
E_L(g,[\xi^{-1},1])$. Hence $E_L(g,[\xi^{-1},1]) = E_L(g\circ \sigma, [-1,1])$.

\par For the sake of simplicity, we let $C := B(1-\xi^{-1})$ and
$\gamma : [-1,1]\rightarrow \R_{+}$ is defined by
$\gamma(x) = \exp\{-C(x+1)\}$. From the discussion in the previous paragraph, we
have indeed reduced the problem to finding $E_L(\gamma, [-1,1])$. This is
because
\begin{equation}
  \label{eq:15}
  E_L(g, [\xi^{-1},1])%
  = E_L(g\circ \sigma, [-1,1])%
  = \exp\{-2B\xi^{-1}\}E_L(\gamma, [-1,1]).
\end{equation}

\par To find a lower bound on $E_L(\gamma,[-1,1])$, we will exploit the
well-known relationship between uniform approximation on the interval by
polynomials and uniform approximation of periodic even functions by
trigonometric polynomials. We write $\conteven[-1,1]$ the space of continuous
and even functions on $[-1,1]$, and for any $L\in \Z_+$ we let
$\tpol_L[-1,1]$ denote the set of even trigonometric polynomials of degree at
most $L$, 
\textit{i.e.} $\tpol_L[-1,1]$ is%
\begin{equation*}
  \left\{ T \in \conteven[-1,1] : \;  T(x) =
    \textstyle\sum_{k=0}^L a_k \cos(\pi k x),\ a_k \in \R,\, x\in [-1,1]\right\} .
\end{equation*}
We furthermore define the periodization operator $P : \cont[-1,1] \rightarrow
\conteven[-1,1]$ such that $Pf( \theta) = f(\cos(\pi \theta))$ for all $f \in
\cont[-1,1]$ and all $\theta \in [-1,1]$. 
Then, it is well-known (see for
instance the Theorem~14.8.1 in \cite{Dav(09)}) that
\begin{equation}
  \label{eq:16}
  E_L(\gamma,[-1,1])%
  = \inf\{\|P\gamma - T\|_{\infty} : \;  T \in \tpol_L[-1,1] \}.
\end{equation}

\par We will now bound the r.h.s. of \eqref{eq:16} by a technique inspired from
\citet{newman1976approximation}, which surprisingly work as well for our
setting. For any $K \in \N$, we define the trigonometric polynomial
$T_K : [-1,1] \rightarrow \bbC$ such that
\begin{equation*}
  T_K(\theta)%
  := e^{i\pi (L+1)\theta}\Big\{\sum_{k=0}^{K-1} e^{i2\pi k
    \theta}\Big\}^2.
\end{equation*}
Then, by orthogonality of the trigonometric polynomials, we have that
\begin{equation}
  \label{eq:17}
  \int_{-1}^{-1}|T_K(\theta)|\,\D \theta%
  = \sum_{j=0}^{K-1}\sum_{k=0}^{K-1} \int_{-1}^1e^{i2\pi(j-k)\theta} \,\D
  \theta%
  = K.
\end{equation}
By definition, for every $\varepsilon > 0$ we can find a $Q \in \tpol_L[-1,1]$
such that $\|P\gamma - Q\|_{\infty} \leq E_L(\gamma, [-1,1]) +
\varepsilon$. Choose such $Q$, and remark that \eqref{eq:17} implies,
\begin{align*}
  \Big|\int_{-1}^1( P\gamma(\theta) - Q(\theta)) T_K(\theta)\, \D \theta
  \Big|%
  &\leq \|P\gamma - Q\|_{\infty} \int_{-1}^1|T_K(\theta)|\, \D \theta\\
  &\leq K\{E_L(\gamma,[-1,1]) + \varepsilon  \}.
\end{align*}
On the other hand remark that $Q$ is a trigonometric polynomial of degree at
most $L$, while $T_K$ is a trigonometric polynomial of degree strictly greater
than $L$. Therefore $Q$ is orthogonal to $T_K$. Moreover, the last display is
true for all $\varepsilon > 0$ and for all $K \in \N$, thus it must be the
case that
\begin{equation}
  \label{eq:18}
  E_L(\gamma, [-1,1])%
  \geq \max_{K\in \N}\frac{1}{K} \Big|\int_{-1}^1P\gamma(\theta)
  T_K(\theta)\, \D \theta \Big|.
\end{equation}

\par Interestingly, we can compute the previous integral. Namely,
\begin{align*}
  \int_{-1}^1P\gamma(\theta) T_K(\theta)\,\D \theta%
  &= \sum_{j=0}^{K-1}\sum_{k=0}^{K-1}\int_{-1}^1 \gamma(\cos (\pi\theta))
    e^{i\pi \theta(L+1 + 2j + 2k)}\,\D \theta\\
  &= \sum_{j=0}^{K-1}\sum_{k=0}^{K-1}\int_0^1 \gamma(\cos(\pi \theta))
    \cos(\pi\theta(L+1+2j + 2k)) \,\D \theta
\end{align*}
The integrals involved in the last
display can be expressed in terms of the modified Bessel function (see \cite[pg. 248]{abramo}) denoted here as
$I_\nu (z)$, which equals 
\[
I_k(z) = \frac{1}{\pi} \int_0^\pi e^{z \cos (t)} \cos (\nu t) \D t
\] 
whenever $\nu =k \in \N$ thanks to \cite[formula 10.32.3]{abramo}.
More precisely, from the above considerations and the fact that the modified Bessel
functions are non--negative, we deduce that
\begin{equation*}
  \Big|\int_{-1}^1P\gamma(\theta)
  T_K(\theta)\, \D \theta \Big|%
  = \sum_{j=0}^{K-1}\sum_{k=0}^{K-1} e^{-C} I_{L+1+2j+2k}(C).
\end{equation*}
\Citet{soni1965inequality} proved that $I_{k+1}(z) \leq I_k(z)$ for all
$k \in \N$ and all $z > 0$. Hence, we obtain from the last display and
\eqref{eq:18} the bound
\begin{equation}
  \label{eq:19}
  E_L(\gamma, [-1,1])%
  \geq \max_{K\in \N} K e^{-C}I_{L+4K}(C).
\end{equation}

In the next lemma, We obtain a bound on the modified Bessel function $z \mapsto I_k(z)$
which remains tighter than the classical bound derived in
\citet{luke1972inequalities} when $z \geq k$. The proof of the lemma is to be
found in Section \ref{sec:proof-pro:4-lower}.
\begin{lem}
  \label{lem:5}
  Assume $k\in \N$ and assume that $z > 8\sqrt{1 + (k/z)^2}$. Then,
  \begin{equation*}
    e^{-z}I_k(z) > \frac{\exp\{ -k^2/(2z)\} }{2e^4(1+(k/z)^2)^{1/4} \sqrt{z}}.%
  \end{equation*}
\end{lem}
\begin{equation*}
  K_{*}%
  :=%
  \begin{cases}
    \alpha \sqrt{C} &\mathrm{if}\ L < \sqrt{C},\\
    \beta C / L &\mathrm{if}\ L \geq \sqrt{C}.
  \end{cases}
\end{equation*}
In view of \eqref{eq:19}, it is clear that $E_L(\gamma, [-1,1]) \geq
K_{*}e^{-C}I_{L+4K_{*}}(C)$. Consider now the case where $L < \sqrt{C}$, then
\begin{equation*}
  0 \leq \frac{L + 4K_{*}}{C}%
  = \frac{L + \alpha \sqrt{C}}{C}%
  < \frac{\alpha + 1}{\sqrt{C}}.
\end{equation*}
Thus, $(L+4K_{*})/C \rightarrow 0$ as $C \rightarrow \infty$, and this implies
that $C > 8\sqrt{1 + (L+4K_{*})^2/C^2}$ when $C$ gets large enough. We then
obtain from Lemma \ref{lem:5} that in this case, as $C \rightarrow \infty$, 
\begin{align*}
  E_L(\gamma, [-1,1])
  &> \frac{\alpha\sqrt{C}(1+o(1))\exp\{-(L+4K_{*})^2/(2C)\}}{2e^2 \sqrt{C}}\\
  &> \frac{\alpha(1+o(1))\exp\{-(\alpha +1)^2/2\}}{2e^2}%
    \gtrsim 1.
\end{align*}

\par We now consider the case $L \geq \sqrt{C}$. In this case, we have,
\begin{equation*}
  0 \leq \frac{L + 4K_{*}}{C}%
  = \frac{L + \beta C/L}{C}%
  = \frac{L}{C} + \frac{\beta}{L}%
  \leq \frac{L}{C} + \frac{\beta}{\sqrt{C}}.
\end{equation*}
Because by assumption there is a constant $\zeta > 0$ such that $L \leq \zeta
C$, then $(L+4K_{*})/C \leq \zeta + o(1)$ as $C \rightarrow \infty$, and thus
we have $C > 8\sqrt{1 + (L+4K_{*})^2/C^2}$ when $C$ is large enough. Then, we
can apply Lemma \ref{lem:5} to find that as $C \rightarrow \infty$, because
$K_{*}^2/C = \beta^2C/L^2 \lesssim 1$ and $K_{*}L/C = \beta \lesssim 1$,
\begin{align*}
  E_L(\gamma,[-1,1])%
  &> \frac{(\beta C /L)(1+o(1))\exp\{-(L + 4K_{*})^2/(2C)\} }{2e^2\sqrt{C}(1 +
    (L/C)^2)^{1/4}}\\
  &= \frac{\beta \sqrt{C}(1+o(1))\exp\{ -L^2/(2C) - 8K_{*}^2/C - 4K_{*}L/C
    \}}{2e^2 L (1+(L/C)^2)^{1/4}}\\
  &\gtrsim \frac{\sqrt{C}(1+o(1))\exp\{ -L^2/(2C) \}}{L
    (1+(L/C)^2)^{1/4}}.
\end{align*}

\par The conclusion then follows from \eqref{eq:15}, by remarking that
$2B\xi^{-1} = 1 + o(1)$ and
$C = B(1-\xi^{-1}) = (\xi/2)(1- \xi^{-1})(1 + O(1/n))$ as
$n \rightarrow \infty$, and thus $L^2/(2C) \leq L^2/\xi + C$ for some $C > 0$ by
definitions of $L$ and $\xi$, and $L/C = 2L(1+o(1))/\xi$.

\subsection{Proof of Lemma \ref{lem:5}}
\label{sec:proof-pro:4-lower}

The proof relies on the well known series representation of the modified Bessel
function (see \cite[formula 10.25.2]{abramo}), namely we have whenever $k \in \N$,
\begin{equation}
  \label{eq:9}
  I_k(z)%
  = \sum_{p=0}^{\infty}\frac{1}{p!(p+k)!}\Big(\frac{z}{2} \Big)^{2p+k}.
\end{equation}
Conveniently, all the terms in the summation are non-negative, which we will
exploit to get our lower bound.

By Stirling's formula, when $k \geq 1$, for any $p \geq 0$
\begin{equation*}
  (p+k)! \leq e\sqrt{(p+k)}\exp\{ -(p+k) + (p+k)\log(p+k)\},
\end{equation*}
and for any $p\geq 1$,
\begin{equation*}
  p! \leq e \sqrt{p}\exp\{-p + p \log p\}.
\end{equation*}
For convenience, let define the functions $\phi_{z,k} : \R_{+} \rightarrow \R_{+}$,
such that for any $x,z \in \R_{+}$ and any $k\in \N$, %
\begin{equation*}
  \phi_{z,k}(x)%
  := -z + 2x + k - x \log x - (x+k)\log(x+k) + (2x + k)\log(z/2).
\end{equation*}
Hence, because each term in the series expansion of \eqref{eq:9} is
non-negative, we get the estimate,
\begin{equation}
  \label{eq:14}
  e^{-z}I_k(z)%
  \geq e^{-z}\sum_{p \geq 1}\frac{1}{p!(p+k)!}\Big(\frac{z}{2} \Big)^{2p+k}%
  \geq \frac{1}{e^2}%
  \sum_{p\geq 1}\frac{\exp\{ \phi_{z,k}(p) \}}{\sqrt{p(p+k)}}.
\end{equation}
Notice that,
\begin{gather*}
  \phi_{z,k}'(x)%
  = - \log(x) - \log(x+k) + 2\log(z/2),\qquad%
  \phi_{z,k}''(x)%
  = -\frac{1}{x} - \frac{1}{x+k}.
\end{gather*}
Thus, $\phi_{z,k}$ admits a unique non-negative extremum at $x_0$ solution to $x_0(x_0
+ k) = z^2/4$, that is,
\begin{equation*}
  x_0%
  = \frac{-k + \sqrt{k^2 + z^2}}{2},\ \mathrm{and},\quad%
  \phi_{z,k}''(x_0) =%
  -\frac{4}{z}\sqrt{1 + (k/z)^2} < 0.
\end{equation*}
Henceforth $x_0$ is indeed the unique maximum of the function $\phi_{z,k}$ on $\R_{+}$. We
let $p_0$ smallest integer larger than $x_0$. Then $p_0 \geq 1$ and we have, by
Taylor expansion that for any $p \geq p_0$ there is a $\bar{p} \in (x_0,p)$
\begin{align*}
  \phi_{z,k}(p)%
  &= \phi_{z,k}(x_0) + \phi_{z,k}'(x_0)(p - x_0) + \frac{1}{2}\phi_{z,k}''(\bar{p})(p -
    x_0)^2\\
  &= \phi_{z,k}(x_0) + \frac{1}{2}\phi_{z,k}''(\bar{p})(p - x_0)^2.
\end{align*}
Remark that, because $\bar{p} \geq x_0$,
\begin{equation*}
  \phi_{z,k}''(\bar{p})%
  = -\frac{1}{\bar{p}} - \frac{1}{\bar{p}+k}%
  \geq  -\frac{1}{x_0} - \frac{1}{x_0+k}%
  = - \frac{4}{z}\sqrt{1 + (k/z)^2}.
\end{equation*}
Then, for any $p \geq p_0$,
\begin{align*}
  \phi_{z,k}(p_0)
  &\geq \phi_{z,k}(x_0) + \frac{1}{2}\phi_{z,k}''(x_0)(p_0 - x_0)^2\\
  &= \phi_{z,k}(x_0) - \frac{2\sqrt{1+ (k/z)^2}}{b}(p - x_0)^2.%
\end{align*}
Therefore,
\begin{equation*}
  e^{-z}I_k(z) \geq \frac{\exp\{\phi_{z,k}(x_0)\}}{e^2}%
  \sum_{p \geq p_0}\frac{\exp\{\phi_{z,k}''(x_0)(p - x_0)^2/2 \} }{\sqrt{p(p+k) }}.
\end{equation*}
Let $p_1$ be the largest integer such that $-\phi_{z,k}''(x_0)(p_1 - x_0)^2 \leq
2$. Remark that whenever $z> 2(1 + (k/z)^2)^{1/2}$, we have $p_1 \geq x_0 + 1$,
which is always the case in the conditions of the lemma. Because the summand is
the previous is monotonically decreasing for $p \geq p_0$, we get the bound, 
\begin{align*}
  e^{-z}I_k(z)
  \geq \frac{\exp\{\phi_{z,k}(x_0)\}}{e^4}\frac{(p_1 -
    p_0)}{\sqrt{p_1(p_1 + k)}}%
  \geq \frac{\exp\{\phi_{z,k}(x_0)\}}{e^4}\frac{(p_1 - x_0) - 1}{\sqrt{p_1(p_1 +
    k)}}.
\end{align*}
But, by the definition of $p_1$, we have that,
\begin{equation*}
  p_1 + 1 - x_0 > \sqrt{\frac{2}{-\phi_{z,k}''(x_0)} }.
\end{equation*}
Therefore, whenever $z > 8(1 + (k/z)^2)^{1/2}$, by the definition of
$\phi_{z,k}''(x_0)$,
\begin{align*}
  e^{-z}I_k(z)%
  &\geq \frac{\exp\{\phi_{z,k}(x_0)\}}{e^4\sqrt{-\phi_{z,k}''(x_0)p_1(p_1 + k)}}%
    \Big\{\sqrt{2} - 2\sqrt{-\phi_{z,k}''(x_0)} \Big\}\\
  &\geq \frac{\sqrt{2}\exp\{\phi_{z,k}(x_0)\}}{2e^4\sqrt{-\phi_{z,k}''(x_0)p_1(p_1 + k)}}.
\end{align*}
Also,
\begin{align*}
  p_1(p_1 + k)
  &= x_0(x_0 + k) + (p_1^2 - x_0^{2}) + (p_1 - x_0)k\\
  &= x_0(x_0 + k) + (p_1 - x_0)(p_1 + x_0 + k)\\
  &=x_0(x_0 + k) + (p_1 - x_0)^2 + (p_1 - x_0)(2x_0 + k).
\end{align*}
But we have that $x_0(x_0 + k) = z^2/4$, $(p_1 - x_0)^2 \leq -2/\phi_{z,k}''(x_0)$,
and $2x_0 + k = z\sqrt{1 + (k/z)^2}$. Thus,
\begin{align*}
  p_1(p_1+k)%
  &\leq \frac{z^2}{4} + \frac{2}{-\phi_{z,k}''(x_0)} +
    \sqrt{\frac{2(1+(k/z)^2)}{-\phi_{z,k}''(x_0)}}z\\
  &= \frac{z^2}{4} + \frac{z}{2\sqrt{1+(k/z)^2}}%
    + \frac{z^{3/2}}{\sqrt{2}}[1 + (k/z)^2]^{1/4}\\
  &= \frac{z^2}{4}\Big\{1 + \frac{z^{-1/2}[1 + (k/z)^2]^{1/4}}{\sqrt{2}}%
    + \frac{z^{-1}}{2\sqrt{1 + (k/z)^2}} \Big\}.
\end{align*}
Therefore, whenever $z > 8(1 + (k/z)^2)^{1/2}$, 
\begin{align*}
  p_1(p_1 + k)%
  \leq \frac{z^2}{4}\Big\{1 + \frac{1}{4} + \frac{1}{16} \Big\}%
  \leq \frac{21}{64} z^2 < \frac{z^2}{2}.
\end{align*}
Hence,
\begin{equation*}
  e^{-z}I_k(z)%
  > \frac{\exp\{\phi_{z,k}(x_0)\} }{e^4\sqrt{-\phi_{z,k}''(x_0)} z}%
  = \frac{\exp\{\phi_{z,k}(x_0)\} }{2e^4(1 + (k/z)^2)^{1/4}  \sqrt{z} }.
\end{equation*}

The remainder of the proof is now dedicated to deriving a lower bound on
$\phi_{z,k}(x_0)$. After some algebra, we find that
\begin{align*}
  \phi_{z,k}(x_0)%
  &= -z + z\sqrt{1 + (k/z)^2}\\
  &\quad
    - (z/2)\{-(k/z) + \sqrt{1+(k/z)^2}\}  \log\{ -(k/z) + \sqrt{1+ (k/z)^2 }
    \}\\
  &\quad - (z/2)\{(k/z) + \sqrt{1 + (k/z)^2}\} \log\{(k/z) + \sqrt{1 + (k/z)^2}
    \}.
\end{align*}
Now we define the function $\varphi : \R_{+} \rightarrow \R$ such that
\begin{multline*}
  \varphi(x)%
  :=%
  -1 + \sqrt{1 + x^2}%
  - \frac{1}{2}(-x + \sqrt{1 + x^2})\log(-x + \sqrt{1+x^2})\\%
  - \frac{1}{2}(x + \sqrt{1 + x^2})\log(x + \sqrt{1+x^2}).
\end{multline*}
Notice that $\phi_{z,k}(x_0) = z \varphi(k/z)$. Also,
\begin{gather*}
  \varphi'(x)%
  = \frac{(- x + \sqrt{1+x^2})\log(-x + \sqrt{1+x^2}) }{2\sqrt{1+x^2}}%
  - \frac{(x + \sqrt{1+x^2})\log(x + \sqrt{1+x^2}) }{2\sqrt{1+x^2}},\\
  \varphi''(x)%
  =-\frac{1}{(1+x^2)^{1/2} },\qquad%
  \varphi'''(x)%
  = \frac{x}{(1 + x^2)^{3/2}}.
\end{gather*}
By a Taylor expansion of $\varphi$ near $0$, we find that there is a $y \in
(0,x)$ such that
\begin{align*}
  \varphi(x)
  &= \varphi(0) + \varphi'(0) x + \frac{1}{2}\varphi''(0) x^2 +
    \frac{1}{6}\varphi'''(y) x^3%
  \geq - \frac{x^2}{2},
\end{align*}
because $\varphi(0) = \varphi'(0) = 0$ and $\varphi'''(y) \geq 0$ for all $y
\geq 0$ by the computations above. This gives the proof for the lower bound on
$\phi_{z,k}(x_0)$ as well, concluding the proof.

\end{document}